\documentclass[11pt]{article}

\usepackage{amsthm}

\newcommand\tsum{\textstyle\sum\nolimits}

\usepackage[nocompress]{cite}
\usepackage{graphicx}

\usepackage{hyperref}
\hypersetup{
    colorlinks=true,
    linkcolor=blue,
    anchorcolor=blue,
    citecolor=blue,
}

\usepackage{enumitem}

 \usepackage{amsmath}
\usepackage{amsfonts}
\usepackage{amssymb, amsfonts, amsmath, nicefrac}
\usepackage{color}
\usepackage{euscript}

\usepackage[margin=1.0in]{geometry}

\def\sV {{\widetilde{V}}}
\def\sW {{\widetilde{W}}}

\newcommand{\cS}{{\mathfrak S}}

\newcommand{\cP}{{\mathfrak P}}

\newcommand{\cM}{{\mathfrak M}}

\newcommand{\cF}{{\mathfrak F}}

\newcommand{\cN}{{\mathfrak N}}

\newcommand{\cR}{{\mathfrak R}}

\newcommand{\cp}{{\mathfrak p}}

\newtheorem{theorem}{Theorem}[section]
\newtheorem{lemma}{Lemma}[section]
\newtheorem{proposition}{Proposition}[section]
\newtheorem{example}{Example}[section]
\newtheorem{remark}{Remark}[section]
\newtheorem{definition}{Definition}[section]

\newtheorem{assumption}{Assumption}[section]
\newtheorem{cor}{Corollary}[section]

\newcommand{\bfR}{{\sf R}}

\def\argmin{\mathop{\rm arg\,min}}
\def\argmax{\mathop{\rm arg\,max}}

\newcommand{\U}{{\cal U}}

\newcommand{\Z}{{\cal Z}}
\newcommand{\F}{{\cal F}}

\newcommand{\A}{{\cal A}}

\newcommand{\G}{{\cal G}}

\newcommand{\X}{{\cal X}}

\newcommand{\PP}{{\cal P}}
\newcommand{\LL}{{\cal L}}

\newcommand{\s}{{\cal S}}
\newcommand{\R}{{\cal R}}
\newcommand{\T}{{\cal T}}

\newcommand{\half}{ \mbox{\small$\frac{1}{2}$}}

\newcommand{\be}{\begin{equation}}
\newcommand{\ee}{\end{equation}}

\def\w{\omega}
\def\e{\varepsilon}
\def\O{\Omega}

\def\supp {{\rm supp}}

\newcommand{\avr}{{\sf AV@R}}

\newcommand{\VaR}{{\sf V@R}}
\newcommand{\RVaR}{{\sf RV@R}}

\def\bbr{{\Bbb{R}}} 
\def\bbe{{\Bbb{E}}} 

\def\bbp{{\Bbb{P}}}

\def\bbb{{\Bbb{B}}}
\def\bbk{{\Bbb{K}}}

\newcommand{\ind}{{\mbox{\boldmath$1$}}}

\begin{document}

\begin{titlepage}
\title{\bf Risk-averse formulations of Stochastic Optimal Control and Markov Decision Processes}
\author{{\bf Alexander Shapiro}\thanks{Georgia Institute of Technology, Atlanta, Georgia
30332, USA, \tt{ashapiro@isye.gatech.edu}\newline
Research of this
author was partially supported by Air Force Office of Scientific Research (AFOSR)
under Grant FA9550-22-1-0244.}
\and
{\bf Yan Li}\thanks{ Texas A\&M University, College Station, Texas 77840, USA, \tt{gzliyan113@tamu.edu}.
}
}


\maketitle
\begin{abstract}
The aim of this paper is to investigate  risk-averse and distributionally robust modeling of Stochastic Optimal Control (SOC) and Markov Decision Process (MDP).  We discuss construction of conditional   nested risk functionals,  a particular attention is given to the Value-at-Risk measure. Necessary and sufficient conditions for existence of non-randomized optimal policies in the framework of robust SOC and MDP are derived. We also investigate sample complexity of optimization problems involving the Value-at-Risk measure.
\end{abstract}


\end{titlepage}

\setcounter{equation}{0}
\section{Introduction}
\label{sec-intr}

The aim of this paper is to investigate  risk averse and distributionally robust approaches to  optimization problems where decisions are made sequentially under conditions of uncertainty. Specifically we discuss  the Stochastic Optimal Control (SOC) and Markov Decision Process (MDP)  formulations  of such optimization problems. For static (single stage) stochastic programs the risk averse and distributionally robust counterparts, of the respective risk neutral problems, are well studied now. Distributionally robust approach to stochastic programming is going back to the pioneering paper by Scarf \cite{scarf}. For a recent survey we can refer to
  \cite{KSW2024},  and   to \cite{esf-2018} where the special attention
 is given   to  construction of the ambiguity sets based on Wasserstein metric.
Modern theory of risk measures was  started in Artzner et al  \cite{ADEH:1999}, where   axioms of the so-called coherent risk measures were formulated. For a discussion of the corresponding risk averse  optimization problems we can  refer to \cite{ruszshap1}. The relation  between distributionally robust and risk averse approaches to stochastic optimization is based on  dual representation of   coherent risk measures (e.g., \cite[Section 6.3]{SDR}).

An extension of risk averse and distributionally robust approaches to sequential optimization is not straightforward and   still  is  somewhat controversial.    An approach to risk averse multistage stochastic programming based on {\em nested} compositions of    coherent risk measures was suggested in \cite{ruzshap2} (see also \cite[Section 6.5.4]{SDR}). That approach involves construction of conditional counterparts of law invariant coherent risk measures. It can be readily extended to SOC problems where the underlying random process does not depend to states and    controls. An extension of the risk averse  to MDPs is more involved since there the evolution of the system is determined  by conditional probabilities (transition kernels) rather than explicitly defined random process. It was suggested by Ruszczy\'nski \cite{rusz2010} to construct nested risk measures directly on the histories of the states  process equipped with the probability law driven  by     (non-randomized) policies of the considered MDP. The approach in \cite{rusz2010} is focused  on coherent risk measures and is     based on their  dual representation.

A different approach to the distributionally robust formulation  of MDPs was suggested in \cite{iyen} and \cite{Nilim2005}, and was used since then  in numerous publications. The approach in   \cite{iyen} and \cite{Nilim2005}  in a sense is static since   the ambiguity sets of transition kernels are defined before realizations of  the decision process. As a consequence,  in order to derive  the corresponding dynamic equations there is a need to introduce the so-called rectangularity conditions on the ambiguity sets of transition kernels. The distributionally robust and risk averse approaches to MDPs can be unified by defining the process as a dynamic  game between the decision maker (the controller) and
the adversary (the nature) (cf., \cite{LiSha}).

The main contribution of this manuscript can be summarized as the following.
\begin{itemize}
     \item
     We discuss construction of conditional and nested risk functionals without assuming their convexity property.   The construction is   general and is not based of the dual (distributionally robust) formulation. A particular attention is given to the Value-at-Risk measure, which is not convex.
   \item We discuss sample complexity of optimization problems involving  the Value-at-Risk measure.
  \item  We give necessary and sufficient conditions for existence of non-randomized optimal policies.
      \item We discuss both the SOC and MDP formulations in  finite and infinite horizons settings. In particular, in the case of MDP our approach is  more general and simpler  than the one of \cite{rusz2010}.
  \end{itemize}

We use the following notation and terminology. For $a\in \bbr$, we denote $[a]_+:=\max\{0,a\}$. For $a,b\in \bbr$ we denote $a\vee b:=\max\{a,b\}$.
A set $\O$ equipped with its sigma algebra $\F$ is called measurable space.
The  measurable space $(\O,\F)$ is said to be {\em  Polish} if $\O$ is  a separable complete metric  space and $\F$ is its Borel sigma algebra. In particular, any   closed subset of   $\bbr^n$ equipped with its Borel sigma algebra  is a Polish space.
 By $\cP$ we denote  the set of probability measures on $(\O,\F)$.
For $P\in \cP$, $(\O,\F,P)$ is called the probability space. It is said  that $P\in \cP$ is nonatomic, if for any $A\in \F$ such that $P(A)>0$, there exists $A'\in \F$ such that $A'\subset A$ and $P(A)>P(A')>0$.
For a set $A\subset \O$  denote by $\ind_A$ its indicator function, i.e., $\ind_A(\w)=1$ for $\w\in A$, and $\ind_A(\w)=0$ for $\w\in \O\setminus A$.
 By $\supp (P)$ we denote the support of probability measure $P$ defined on a metric space $\O$ equipped with its Borel sigma algebra, i.e.,
$A=\supp (P)$ is the smallest closed subset of $\O$ such that $P(A)=1$.

For $P\in \cP$ and  a random variable    $Z:\O\to\bbr$   we denote by $\bbe^P[Z]=\int_\O ZdP$ its expected value, and by
$$F_Z^P(z):=P(Z\le z),\;z\in \bbr,$$ its cumulative distribution function (cdf). By  $\cF$ we denote the set of right side continuous monotonically nondecreasing functions $\phi:\bbr\to \bbr$ such that $\lim_{z\to-\infty}\phi(z)=0$ and $\lim_{z\to +\infty}\phi(z)=1$, i.e. $\cF$ is the space of cdfs.
By $L_p(\O,\F,P)$ we denote the space of random variables $Z$ such that $\int_\O|Z|^p dP<\infty$, $p\in [1,\infty)$.

By $\delta_\xi$ we denote measure of mass one at a point $\xi$ (Dirac measure).
For $Q,P\in \cP$ it is said that $Q$ is absolutely continuous with respect to $P$, denoted $Q\ll P$, if for $A\in \F$ such that  $P(A)=0$  it follows that $Q(A)=0$.
For a sigma subalgebra $\G$ of $\F$ we denote by $\bbe^P_{|\G}[Z]$  the respective  conditional expectation. That is,  $\bbe^P_{|\G}[Z]$ is $\G$-measurable and $P$-integrable, and\footnote{This is the classical definition of conditional expectation due to Kolmogorov.}
\[
\int_A \bbe^P_{|\G}[Z] dP=\int_AZdP, \;\forall A\in \G.
\]
Note that $\bbe^P_{|\G}[Z]$ consists of a family of $\G$-measurable random variables (called versions) which are equal to each other $P$-almost surely.
By $P_{|\G}(A)=\bbe_{|\G}[\ind_A]$, $A\in \F$, is denoted the respective conditional probability.
For a random variable $Y$, we denote by $\bbe_{|Y}[Z]$  the conditional expectation and by $P_{|Y}(A)=\bbe_{|Y}[\ind_A]$, $A\in \F$,   the respective conditional probability.

\setcounter{equation}{0}
\section{Risk functionals}
\label{sec-riskfun}

Let $(\O,\F)$ be a measurable space. Assume that with every probability measure  $P\in \cP$  is associated a linear space $\Z$  of measurable functions (random variables) $Z:\O\to \bbr$ and functional $\R:\Z\to \bbr$.  Consider the following axioms (conditions) that $\R$  may  satisfy.
\begin{itemize}
  \item [(A1)] (monotonicity) if $Z,Z'\in \Z$ and $Z\ge Z'$, $P$-a.s., then $\R(Z)\ge \R(Z')$.
    \item [(A2)] (convexity)    if $Z,Z'\in \Z$ and $\tau\in [0,1]$, then
    $\R(\tau Z+(1-\tau)Z')\le \tau \R(Z)+(1-\tau)\R(Z')$.
    \item [(A3)] (translation equivariance) for   $Z\in \Z$ and $\tau\in \bbr$, it follows that $\R(Z+\tau)=\R(Z)+\tau$.
     \item [(A4)] (positive homogeneity) if $Z\in \Z$ and $\tau\ge 0$, then
     $\R(\tau Z)=\tau \R(Z)$.
      \end{itemize}
In order to emphasize that the functional $\R$ depends on $P\in \cP$,  we write it as $\R^P$. The corresponding linear space $\Z$ may also depend on $P$, but we suppress this in the notation. In any case we assume that if $Z\in \Z$ and $\tau\in \bbr$, then $Z+\tau\in \Z$.
Note that it follows from axiom (A1) that if $Z=Z'$, $P$-a.s., then $\R^P(Z)= \R^P(Z')$.
It follows from axiom (A4) that $\R^P(0)=0$.

Functionals satisfying axioms (A1) - (A4) are called {\em coherent} risk measures, \cite{ADEH:1999}.
An important example of coherent risk measure is the Average Value-at-Risk (also called Conditional  Value-at-Risk, expected shortfall,   expected tail loss). The following representation of the Average Value-at-Risk is due to \cite{ury2},
\begin{equation}\label{avar}
\avr^P_{\alpha}(Z):=\inf_{\tau\in \bbr}\{\tau+\alpha^{-1}\bbe^P[Z-\tau]_+\},\;\alpha\in (0,1].
\end{equation}
In that example, the respective  space $\Z=L_1(\O,\F,P)$.
An important example of convex functional, satisfying axioms (A1) -  (A3), is the entropic risk measure
 \[
 \R^P_{\tau}(Z):=\tau^{-1}\log \bbe^P[e^{\tau Z}], \;\tau>0.
 \]
 Here the respective space $\Z$ consists of random variables $Z$ such that $\bbe^P[e^{\tau Z}]<\infty$ for all $\tau\in \bbr$.

 An important example of non-convex risk measure   is the Value-at-Risk, $\R^P=\VaR^P_{\alpha}$, where
\begin{equation}\label{quan}
\VaR^P_{\alpha}(Z):=\inf\{z:F_Z^P(z)\ge 1-\alpha\},\;\;\alpha\in (0,1).
\end{equation}
That is, $\VaR^P_{\alpha}(Z)$ is the left side $(1-\alpha)$-quantile of the distribution of $Z$. In that example the corresponding space $\Z$ consists of all measurable $Z:\O\to \bbr$.
The $\VaR^P_{\alpha}(\cdot)$ functional satisfies axioms (A1),(A3) and (A4), but not axiom (A2).

\begin{itemize}
  \item []
  It is said that $\R^P$ is {\em strictly monotone} if $Z\ge Z'$, $P$-a.s.,  and $P(Z>Z')>0$, then $\R^P(Z) > \R^P(Z')$.
The $\avr^P_{\alpha}$ and  $\VaR^P_{\alpha}$, $\alpha\in (0,1)$,  functionals are monotone but are  not strictly monotone.
\end{itemize}

 It is said that $Z,Z'\in \Z$ are distributionally equivalent (with respect to $P\in \cP$), if their cdfs do coincide, i.e.,
 $P(Z\le z)=P(Z'\le z)$ for all $z\in \bbr$.
\begin{definition}
It is said that $\R^P$ is {\rm law invariant}   if for any distributionally equivalent  $Z,Z'\in \Z$ it follows that  $\R^P(Z)=\R^P(Z')$.
\end{definition}

That is,  a law invariant risk measure $\R^P(Z)$ is a function of the cdf $F_Z^P$, i.e.,
can be represented as
 \begin{equation}\label{cdfr}
   \R^P(Z)=\rho (F_Z^{P}),\;\;Z\in \Z,
 \end{equation}
where $\rho$ is mapping the corresponding $\phi\in \cF$ into $\bbr$.
For example, for $\R^P=\avr^P_\alpha$ the corresponding
\[
\rho(\phi)=\inf_{\tau\in \bbr}\left\{\tau+\alpha^{-1} \int_{-\infty}^{+\infty} [z-\tau]_+d \phi(z)\right\},
\]
defined for such $\phi \in \cF$ that $\int_{-\infty}^{+\infty} [z-\tau]_+d \phi(z)<\infty$ for all $\tau\in \bbr$.
For $\R^P=\VaR_\alpha^P$ the corresponding mapping is
$$\rho(\phi)=\inf\{z:\phi(z)\ge 1-\alpha\},\;
\phi\in \cF.$$


We also consider robust counterparts of risk functionals. That is, let $\cM\subset \cP$  be a (nonempty)  set of probability measures. Then the robust counterpart of $\R^P$ is defined as
\begin{equation}\label{robcoun}
  \bfR(Z):=\sup_{P\in \cM} \R^P(Z),
\end{equation}
assuming that $\sup_{P\in \cM} \R^P(Z)<\infty$ for all $Z\in \Z$.
If $\R^P$ satisfies any of axioms (A2) - (A4) for all $P\in \cM$, then so is its
robust counterpart $\bfR$. The monotonicity axiom is more involved, this is because the inequality $Z\ge Z'$ is defined $P$-almost surely  with respect to a particular
$P\in \cP$. Note  that  if $Z\ge Z'$ almost surely with respect to $\bbp$, then $Z\ge Z'$ almost surely with respect to any measure  $P\in \cP$    absolutely continuous with respect to $ \bbp$. Therefore we   consider the following setting. Consider a probability measure $\bbp\in \cP$, viewed as a reference measure. Suppose that $\bbp\in \cM$,  and every $P\in \cM$ is absolutely continuous with respect to $\bbp$.
Suppose that for every $P\in \cM$ the functional $\R^P$ is monotone. Then the robust functional is monotone with respect to $\bbp$, i..e, if $Z\ge Z'$, $\bbp$-a.s., then $\bfR(Z)\ge \bfR(Z')$.

In particular   let $\R^P:=\bbe^P$.    Then the corresponding    functional $\bfR$, referred to as
  the  {\em distributionally robust functional}, becomes
\begin{equation}\label{robustfun}
   \bfR(Z)=\sup_{P\in \cM}\bbe^P(Z).
\end{equation}
Of course, $\bbe^P$ is monotone  for every $P\in \cP$, and hence the distributionally robust functional
$\bfR$ is monotone with respect to $\bbp$.
Note that if  the distributionally robust functional $\bfR$  is    law invariant  with respect to\footnote{That is,  if $\bbp(Z\le z)=\bbp(Z'\le z)$ for all $z\in \bbr$, then $\bfR(Z)=\bfR(Z')$.}
$\bbp$, then   every $P\in \cM$  is absolutely continuous with respect to  $\bbp$.  Indeed if $A\in \F$ is such that $\bbp(A)=0$, then $\ind_A$ is distributionally equivalent (with respect to $\bbp)$  to 0, and hence by the law invariance it follows that  $\bfR(\ind_A)=0$.  Consequently  $P(A)\le \bfR(\ind_A)=0$ for any $P\in \cM$.

\begin{example}[Robust Value-at-Risk]
Consider $\VaR^P_{\alpha}$. The  the corresponding  robust  Value-at-Risk is
\begin{equation}\label{robvar}
\begin{array}{lll}
 \RVaR_{\alpha}(Z)&:=&\sup_{P\in \cM}\VaR^P_{\alpha}(Z)=
  \inf\left \{z:P(Z\le z)\ge 1-\alpha,\;\forall P\in \cM\right \}\\
  &=&  \inf\left \{z:\inf_{P\in \cM}P(Z\le z)\ge 1-\alpha\right \}.
 \end{array}
\end{equation}
Consider a {\em nonatomic}  probability measure $\bbp\in \cP$, viewed as a reference measure.
Suppose that $\bbp\in \cM$ and  the respective distributionally robust functional $\bfR(Z)$, defined in \eqref{robustfun},   is {\em law invariant} with respect to $\bbp$.
 We have  that  for sets $A,A'\in \F$,   the indicator functions $Z=\ind_A$ and $Z'=\ind_{A'}$  are distributionally equivalent with respect to $\bbp$  iff $\bbp(A)=\bbp(A')$.  Therefore  if $\bbp(A)=\bbp(A')$, then
$
  \bfR(\ind_A)=\bfR(\ind_{A'}).
$
 Consider function $\cp:[0,1] \to \bbr$ defined as
 \begin{equation}\label{funcp}
  \cp(\tau):= \inf_{P\in \cM} P(A)\;\text{ for some}\;A\in \F\;\text{such that } \;\bbp(A)=\tau.
 \end{equation}
 Since $P(A)=\bbe^P[\ind_A]$ and because of the assumed law invariance of  the respective  distributionally robust functional,  $\cp(\tau)$ is well defined; and since $\bbp$ is nonatomic, $\cp(\tau)$ is defined for every $\tau\in [0,1]$.

 Properties of function $\cp(\tau)$  are discussed in
\cite[Proposition 7.16]{SDR}. In particular, $\cp(\tau)\le \tau$, and   $\cp(\cdot)$ is monotonically nondecreasing  and continuous on the interval $(0,1]$. We have then (cf., \cite[Prpoposition 7.14]{SDR})
\begin{equation}\label{dynsoc-4}
\RVaR_{\alpha}(Z)=  \inf\left \{z: \bbp(Z\le z)\ge \cp^{-1}(1-\alpha)\right \}.
\end{equation}
That is, $\RVaR_{\alpha}(\cdot)=\VaR^\bbp_{\alpha^*}(\cdot)$, where $1-\alpha^*=\cp^{-1}(1-\alpha)$.
Note that  $1-\alpha^*\ge 1- \alpha$.
\end{example}

\subsection{Sample complexity }\label{sec_sample_var}

Let $Z_1,...,Z_N$  be iid samples of random variable $Z\sim P$,  and $P_N$  be the corresponding empirical measure, i.e., $P_N=N^{-1}\sum_{i=1}^N \delta_{Z_i}$. Furthermore, let   $F^{P_N}(z)=P_N(Z\le z)=N^{-1}\sum_{i=1}^N\ind_{\{Z_i\le z\}}$ be the empirical cdf and
\[
\VaR^{P_N}_{\alpha} :=\inf\left\{z:F^{P_N}(z)\ge 1-\alpha\right \},
\]
be the    empirical Value-at-Risk viewed as an estimator of   $\VaR^{P}_{\alpha}(Z)$.
We are interested in evaluating the sample size $N$ such that $\left|\VaR^{P_N}_{\alpha}-\VaR^{P}_{\alpha}(Z)\right |<\e$ with   high probability  for a given $\e>0$.
Note that if the equation $F^P_Z(z)=1-\alpha$ has more than one solution, i.e., the $(1-\alpha)$-left side quantile  is smaller than the  $(1-\alpha)$-right side quantile, then $\VaR^{P_N}_{\alpha}$ may not converge,  in probability,  to $\VaR^{P}_{\alpha}(Z)$    as $N\to\infty$. That is, in general there is no guarantee that $\VaR^{P_N}_{\alpha}$ is a consistent estimator of $\VaR^{P}_{\alpha}(Z)$.

Recall the Dvoretzky-Kiefer-Wolfowitz (DKW) inequality  (e.g., \cite[Theorem 11.6]{Kosorok}): for $Z\sim P$,
$F^P =F_Z^P$ and $\e>0$  the following inequality holds
\begin{equation}\label{dkf}
P\left(\sup_{z \in \bbr} \left|{F^{P_N}(z) - F^{P}(z) }\right| >\e\right)\leq 2 e^{-2 N \e^2}.
\end{equation}
Consider first the case where $P\left (Z=\nu\right )>0$, with  $\nu:=\VaR^{P}_{\alpha}(Z)$,  i.e.\footnote{Note that  $F^P(\cdot)$ is right side continuous and $\lim_{z\downarrow \nu}F^P(z)=F(\nu)$.},
 $\lim_{z\uparrow \nu}F^P(z)<\lim_{z\downarrow \nu}F^P(z) .$
Moreover suppose that
\begin{equation}\label{epsilon}
 \begin{array}{l}
\kappa:=\min\left\{ (1-\alpha)-\lim_{z\uparrow \nu}F^P(z), \lim_{z\downarrow \nu}F^P(z)-(1-\alpha)\right\}
>0.
\end{array}
\end{equation}
Then by the DKW inequality we have that
\begin{equation}\label{sample-1}
  P\left(\VaR^{P_N}_{\alpha} =\VaR^{P}_{\alpha}(Z)\right )\ge 1-2e^{-2 N \kappa^2}.
\end{equation}
This implies the following.
\begin{proposition}\label{prop_cvar_saa_exact}
Suppose that condition \eqref{epsilon}  holds. Then for $\delta>0$ and $N\ge \half \kappa^{-2}\log(2/\delta)$ the empirical  $\VaR^{P_N}_{\alpha}$ is {\em equal} to
$\VaR^{P}_{\alpha}(Z)$  with probability at least $1-\delta$.
\end{proposition}

Another case which we consider is the following. Suppose that $F^P=F^P_Z$ has the following local  growth property.
\begin{itemize}
  \item
  There are constants $b>0$ and $c>0$ such that
\begin{equation}\label{growth}
 F^P(z')-F^P(z)\ge  c\, (z'-z) \;\text{for all}\; z,z' \in [\nu-b,\nu+b]\;\text{such that}\; \;z'\ge z,
\end{equation}
 where  $\nu:=\VaR^{P}_{\alpha}(Z)$.
\end{itemize}
 In particular this condition  holds if the density $dF^{P}(z)/dz$ exists  and $dF^{P}(z)/dz\ge c$ for all $z\in [\nu-b,\nu+b]$.

 For $\phi\in \cF$ and $\alpha\in (0,1)$ denote $$\VaR_{\alpha,\phi}:=\inf\{z:\phi(z)\ge 1-\alpha\}.$$
 For $\epsilon>0$ and $\phi\in \cF$
suppose that   $$\sup_{z\in [\nu-b,\nu+b]}\left|\phi(z)  - F^{P}(z) \right| \le \epsilon.$$  It follows by the growth condition \eqref{growth}  that if $\epsilon <b c$, then  $\VaR_{\alpha,\phi}\in [\nu-b,\nu+b]$ and
\begin{equation}\label{growth-2}
\left|\VaR^{P}_{\alpha}(Z)-\VaR_{\alpha,\phi}\right |\le \epsilon/c,
\end{equation}
By the DKW inequality this implies the following.

\begin{proposition}
\label{prop_concentration_var}
Suppose that the growth condition \eqref{growth} holds and let $\delta>0$.
Then   for any $\e\in (0,b)$ and  $N\ge \half  c^{-2}  \e^{-2}  \log(2/\delta)$  it follows
\begin{equation}\label{sample-3}
 P\left(|\VaR^{P_N}_{\alpha} (Z)  -\VaR^{P}_{\alpha}(Z)|<\e\right )\ge 1-\delta.
\end{equation}
\end{proposition}

  \begin{proof}
  We have that $\VaR^{P_N}_{\alpha} \in [\nu-b,\nu+b]$ if
  \[
  \sup_{z\in [\nu-b,\nu+b]}\left|F^{P_N}(z)  - F^{P}(z) \right| \le bc.
  \]
  By the DKW inequality and \eqref{growth-2},
  for  $N \ge \half  c^{-2}  b^{-2} \log(2/\delta)$,
\begin{align*}
P\left(\left|\VaR^{P_N}_{\alpha}(Z)-\VaR^{P}_{\alpha}(Z)\right|  < \min \left\{b , \sqrt{ \tfrac{\log(2/\delta))}{2N c^2}}  \right\} \right ) \geq  1-\delta.
\end{align*}
  Thus  we have  that
  for
   \begin{equation}\label{sample-5}
   N\ge\left (\half  c^{-2}  b^{-2} \log(2/\delta)\right)\vee \left (\half  c^{-2}  \e^{-2}  \log(2/\delta)\right )
   \end{equation}
  the bound \eqref{sample-3} follows.
Of course for $\e<b$ the first term in the right hand side of \eqref{sample-5} can be omitted, and hence the proof is complete.
\end{proof}

Consider now  the following optimization problem
\begin{equation}
\label{probvar}
\min_{x \in {\cal X}} \VaR_\alpha^{P} \left (\psi(x, \xi) \right),
\end{equation}
where $\X$ is a (nonempty)  compact subset of $\bbr^n$,
$\xi$ is a random vector whose probability distribution $P$ is supported on a closed set  $\Xi\subset \bbr^d$,  and
$\psi:\X\times \Xi\to \bbr$. Denote   $Z_x (\cdot) := \psi(x, \cdot)$.
We assume that $Z_x:\Xi\to \bbr$ is  (Borel) measurable for every $x\in \X$.
 We say that the  growth condition \eqref{growth} holds {\em uniformly} if   there are constants $b>0$ and $c>0$ such that for every $x\in \X$,
 \begin{equation}\label{growth-2a}
 F_x(z')-F_x(z)\ge  c\, (z'-z) \;\text{for all}\; z,z' \in [\nu_x-b,\nu_x+b]\;\text{such that}\; \;z'\ge z,
\end{equation}
 where $F_x$ is the cdf of $Z_x$,  and
 $\nu_x:=\VaR^{P}_{\alpha}(Z_x)$.
 \begin{assumption}
\label{ass-lip}
There is  a positive constant $L$ such that
\begin{equation}
\label{condition_lipscthiz_in_x}
 \left| \psi(x, \xi) - \psi(x', \xi) \right| \leq L \left\Vert x - x' \right\Vert, \;\forall
 (x,\xi)\in \X\times \Xi.
\end{equation}
\end{assumption}

Let $\xi_1,...,\xi_N$  be an iid sample of $\xi\sim P$,  and $P_N$ be the corresponding empirical measure.
Denote $D:=\sup_{x,x'\in \X}\|x-x'\|$ the diameter of the set $\X$. Note that $D$ is finite, since $\X$ is assumed to be compact and hence is bounded.

\begin{theorem}
\label{pr-unifbound}
Suppose Assumption {\rm \ref{ass-lip}} and the uniform growth condition \eqref{growth-2a} hold.
Then for any $\e  \in (0,b)$, $\delta>0$  and
\begin{equation}\label{sampeps-1}
 N  \geq 2  c^{-2} \e^{-2} \left[ n \log \left( \tfrac{ 4 L D}{ \e  } \right) + \log(\tfrac{1}{\delta}) \right],
\end{equation}
it follows   that
\begin{equation}\label{sampeps-2}
P\left( \sup_{x\in \X}\left| \VaR_\alpha^{P} \left (\psi(x, \xi) \right) - \VaR_\alpha^{P_N} \left (\psi(x, \xi) \right) \right|
 \leq \e\right ) \ge 1-\delta.
\end{equation}
 \end{theorem}

\begin{proof}
We begin by first showing that for $Z_x (\cdot) = f(x, \cdot)$,
\begin{align*}
\left| \VaR_\alpha^{P} \left (Z_x \right) - \VaR_\alpha^{P} \left(Z_{x'} \right) \right| \leq L \left\Vert x - x' \right\Vert,\;x,x'\in \X.
\end{align*}
Indeed, let $\nu_\alpha^x := \VaR_\alpha^{P} (Z_x)$. It is clear that for any $\epsilon > 0$, we have
$
P(Z_x \leq \nu_{\alpha}^x + \epsilon) \geq 1-\alpha,
$
and hence
\begin{align*}
P(Z_{x'} - \left\Vert Z_{x'} - Z_x \right\Vert \leq\nu_{\alpha}^x + \epsilon) \geq 1-\alpha.
\end{align*}
This implies $\nu_\alpha^{x'} \leq \nu_\alpha^x + \left\Vert Z_{x'} - Z_x \right\Vert + \epsilon$. Since this holds for any $\epsilon >0$, it follows that
$\nu_\alpha^{x'} \leq \nu_\alpha^x + \left\Vert Z_{x'} - Z_x \right\Vert$.
Interchanging $\nu_\alpha^{x'}$ and $\nu_\alpha^x$ yields
\begin{align}\label{condition_lipscthiz_in_x_sample-1}
\left| \nu_\alpha^{x'} - \nu_\alpha^x \right| \leq \left\Vert Z_{x'} - Z_x \right\Vert \leq L \left\Vert x - x' \right\Vert,\;x,x'\in \X,
\end{align}
where the last inequality follows from \eqref{condition_lipscthiz_in_x}.
Note the same treatment also shows
\begin{align}\label{condition_lipscthiz_in_x_sample}
\left| \VaR_\alpha^{P_N} \left (Z_x \right) - \VaR_\alpha^{P_N} \left(Z_{x'}) \right) \right| \leq L \left\Vert x - x' \right\Vert,\;x,x'\in \X.
\end{align}

Suppose the uniform  growth condition \eqref{growth-2a} holds. Let $\eta:=\e/(4L)$ and
consider an $\eta$-net of ${\cal X}$, denoted by ${\cN}_{\eta}$. That is,
for any $x \in {\cal X}$, there is  $x_\eta \in {\cN}_{\eta}$ such that $\left\Vert x - x_\eta\right\Vert \leq \eta$.
For $x\in \X$ we can write
\begin{equation}
\label{ineq}
\begin{array}{lr}
&\left| \VaR_\alpha^{P} \left (Z_x \right) - \VaR_\alpha^{P_N} \left (Z_x \right) \right|
 \leq
\left| \VaR_\alpha^{P} \left (Z_{x_\eta}) \right) - \VaR_\alpha^{P_N} \left (Z_{x_\eta}  \right) \right| +  \\
 & \left |\VaR_\alpha^{P} \left (Z_x \right) - \VaR_\alpha^{P} \left (Z_{x_\eta} \right) \right|
  +
\left| \VaR_\alpha^{P_N} \left (Z_x   \right) - \VaR_\alpha^{P_N} \left ( Z_{x_\eta}  \right) \right|.
\end{array}
\end{equation}
By \eqref{condition_lipscthiz_in_x_sample-1} and \eqref{condition_lipscthiz_in_x_sample} we have that the sum of the last two terms in \eqref{ineq}  is $\le \e/2$.
It follows that
\begin{equation}
\label{probound}
P\left( \sup_{x\in \X}\left| \VaR_\alpha^{P} \left (Z_x \right) - \VaR_\alpha^{P_N} \left (Z_x \right) \right|
 \leq \e\right )\ge P\left(\max_{x_\eta\in \cN_\eta}\left| \VaR_\alpha^{P} \left ( Z_{x_\eta} \right) - \VaR_\alpha^{P_N} \left (Z_{x_\eta} \right) \right|\le \e/2\right).
\end{equation}

Now by Proposition \ref{prop_concentration_var}, for every $x_\eta\in \cN_\eta$  we have that
\[
\left| \VaR_\alpha^{P} \left (Z_{x_\eta} \right) - \VaR_\alpha^{P_N} \left (Z_{x_\eta} \right) \right|<\e/2
\]
with probability at least $1-\delta$ for sample size $N\ge 2  c^{-2}  \e^{-2}   \log(2/\delta)$. Consequently
\[
P\left( \max_{x_\eta\in \cN_\eta}
\left| \VaR_\alpha^{P} \left (Z_{x_\eta} \right) - \VaR_\alpha^{P_N} \left (Z_{x_\eta} \right) \right| \le \e/2 \right)\ge 1- M \delta,
\]
 where   $M:=|\cN_\eta|$ is the cardinality of the net.
By \eqref{probound} it follows that for   $N\ge 2  c^{-2}  \e^{-2}   \log\left(\tfrac{M}{2\delta}\right)$, the bound \eqref{sampeps-2} follows.

We have that cardinality $M=|\cN_\eta|$ of the net   can be bounded by  $C (D/\eta)^n =C(4LD/\e)^n$  for some constant $C>0$. Choice of $C$ may depend on the considered norm $\|\cdot\|$. For example for the $\ell_\infty$ norm (the max-norm) we can use $C=1$. Exact value of the constant $C$ is not important here. We can use generically $C=2$. We obtain that \eqref{sampeps-2} holds for the sample size
$N\ge 2  c^{-2}  \e^{-2}   \log\left(\tfrac{(4 LD/\e)^n}{\delta}\right)$. This completes the proof.
\end{proof}

Denote by $\s_\e$ the set of $\e$-optimal solutions of problem  \eqref{probvar}, i.e.,
\[
\s_\e=\left\{x\in \X: \VaR_\alpha^{P} \left (\psi(x, \xi)\right)\le \inf_{x\in \X}\VaR_\alpha^{P} \left (\psi(x, \xi)\right)+\e\right\},
\]
and let $\hat{\s}_N:=\argmin_{x\in \X}\left\{\VaR_\alpha^{P_N} \left (\psi(x, \xi) \right)\right\}$  be the set of optimal solutions of the respective empirical (Sample Average Approximation) problem. By the uniform bound of Theorem   \ref{pr-unifbound} we have the following (this can be compared with \cite[Theorem 5.18]{SDR}).

\begin{cor}
Suppose Assumption \ref{ass-lip} and the uniform growth condition \eqref{growth-2a} hold.
Then for any $\e  \in (0,b)$, $\delta>0$  and $N$ satisfying \eqref{sampeps-1}, it follows that  the event $\{\hat{\s}_N\subset \s_\e\}$ happens with probability at least   $1-\delta$.
\end{cor}

\subsection{Conditional   risk  functionals}
\label{sec-cond}
We discuss in this section conditional counterparts of risk functionals.
For random variable  $Z :\O\to\bbr$  and sigma subalgebra $\G\subset \F$  consider the conditional cdf
\[
F^P_{Z|\G}(z):=P_{|\G}(Z\le z)=\bbe_{|\G}^P[\ind_{\{Z\le z\}}],\;z\in \bbr.
\]
Suppose that $\R^P$ is law invariant and hence can be represented in the form \eqref{cdfr}. Then the conditional counterpart of $\R^P$ is defined as
\begin{equation}\label{condrisk}
\R^P_{|\G}(Z):=\rho(F^P_{Z|\G}).
\end{equation}
For law invariant coherent risk measures such definition of the respective conditional counterparts is equivalent to the standard definitions used in the literature. For example, in the definition \eqref{avar} of the Average Value-at-Risk  the minimum is attained at $\bar{\tau}=\VaR_\alpha^P(Z)$  and hence
its conditional counterpart can be obtained by replacing $\VaR_\alpha^P(Z)$ with its conditional counterpart.  For the Value-at-Risk its conditional counterpart is
\begin{equation}\label{quan-2}
\VaR^P_{\alpha|\G}(Z)=\inf\left\{z:P_{|\G}(Z\le z)\ge 1-\alpha\right\},
\;\;\alpha\in (0,1).
\end{equation}

In particular we consider the following   construction   where we follow   \cite[Appendix A1 - A2]{ShaPich2024}.
Assume that the measurable space $(\O,\F)$  is given as the product of   measurable spaces $(\O_1,\F_1)$ and $(\O_2,\F_2)$, i.e., $\O=\O_1\times \O_2$ and $\F=\F_1\otimes \F_2$. Denote by $\cP_1$ and $\cP_2$ the  sets of probability measures on the respective measurable spaces  $(\O_1,\F_1)$ and $(\O_2,\F_2)$.
For a probability measure $P\in \cP$ we denote by $P_1\in \cP_1$   the respective marginal probability measure on $(\O_1,\F_1)$, that is
$P_1(A)=P(A\times \O_2)$ for $A\in \F_1$. The marginal probability measure $P_2\in \cP_2$   is defined in the similar way.
We assume that $(\O_1,\F_1)$ and  $(\O_2,\F_2)$  are {\em Polish} spaces.

Consider sigma subalgebra $\G$ of $\F$  consisting  of the   sets $A\times \O_2$, $A\in \F_1$, that is
\begin{equation}\label{subal}
	\G:=\{A\times \O_2: A\in \F_1\}.
\end{equation}
Note that the elements of  subalgebra $\G$  are determined by sets (events) $A\in \F_1$, and in that sense $\G$ can be identified with the sigma algebra $\F_1$, we write this as $\G\equiv \F_1$.
A random variable $Z(\w_{1},\w_{2})$, $\w=(\w_{1},\w_{2})\in \O_1\times\O_2$, is $\G$-measurable iff $Z(\w_{1},\cdot)$ is constant on $\O_2$  for every $\w_{1}\in \O_1$, and $Z(\cdot,\w_{2})$ is $\F_1$-measurable   for every $\w_{2}\in \O_2$.
Therefore, with some abuse of the notation, we write a $\G$-measurable variable as a function $Z(\w_{1})$ of $\w_{1}\in \O_1$. We also use notation  $Z_{\w_{1}}(\w_{2}):=Z(\w_{1},\w_{2})$ viewed as random variable
$Z_{\w_{1}}:\O_2\to \bbr$.


\begin{definition}[Regular Probability Kernel]\label{def-ker}
	A function $\bbk:\F_2\times \O_1\to [0,1]$ is said to be a \emph{Regular Probability Kernel} (RPK) of a probability measure $P\in \cP$ if the following properties hold:
		{\rm (i)} $\bbk(\cdot |\w_{1})$ is a probability measure for $P_1$-almost every $\w_{1}\in \O_1$,
		{\rm (ii)}  for every $B\in  \F_2$ the function $\bbk(B|\cdot)$ is $\F_1$-measurable,
		{\rm (iii)} for every $A\in \F_1$ and $B\in \F_2$ it follows that
	\begin{equation}\label{regker-1}
	  P(A\times B)=\int_A \bbk(B|\w_{1})d P_1 (\w_{1}).
	\end{equation}
	In particular, $P_2(B)=\int_{\O_1} \bbk(B|\w_{1})d P_1 (\w_{1})$ is the respective marginal probability measure.
\end{definition}

The Disintegration Theorem (e.g.,   \cite[III-70]{Dellacherie1978}) ensures existence of the RPK for a wide class of measurable spaces, in particular if
$(\O_1,\F_1)$ and  $(\O_2,\F_2)$   are  Polish spaces.
Therefore, we assume existence of the RPK for every $P\in \cP$.
The function  $\bbk(B|\cdot)$ is defined up to $P_1$-measure zero, and   is uniquely determined on the      $\supp(P_1)$.
The RPK is associated with a specified $P\in \cP$,  we sometimes write $\bbk^P$    to emphasize this.

For $P\in \cP$ and $\w_{1}\in \supp (P_1)$   we can  define a probability measure $P_{|\w_{1}}\in \cP_2$ as
\begin{equation}\label{rpk-1}
P_{|\w_{1}}(B):=\bbk^P(B|\w_{1}),\;B\in \F_2.
\end{equation}
Then
we can define the conditional 
 counterpart of the
cdf $F_Z^P$ as
\begin{equation}\label{condcdf}
 F^P_{Z|\w_1}(z):=P_{|\w_{1}}(Z_{\w_1}\le z), \;P_1- {\rm  a.e.} \;\;\w_1\in \O_1.
\end{equation}
That is, $ F^P_{Z|\w_1}$  is the cdf of $Z_{\w_1}:\O_2\to \bbr$
 associated with the
probability measure $P_{|\w_{1}}$.
The corresponding  conditionally law
invariant risk functional   is   a function of
the conditional cdf. That is, the conditional law
invariant risk functional
$\R^P_{|\w_1}:\Z\to \bbr$ is defined as
\begin{equation}\label{condcdf-2}
\R^P_{|\w_1}(Z)=\rho_1\big (F^P_{Z|\w_1}\big),
\end{equation}
where  $\rho_1$ is the corresponding mapping.
Note that for sigma algebra $\G\equiv \F_1$ defined in \eqref{subal},   $F^P_{Z|\w_1}$ can be viewed as   version   $F^P_{Z|\w_1}=F^P_{Z|\G}(\w_1)$ of the conditional cdf, and
 $\R^P_{|\w_1}$ can be viewed as   version $\R^P_{|\w_1}=\R^P_{|\G}(\w_1)$ of the  conditional risk measure.
For example, for the Value-at-Risk its conditional counterpart
\begin{equation}\label{condcdf-3}
\VaR^P_{\alpha|\w_1}(Z)=\inf\left\{z:F^P_{Z|\w_1}(Z\le z)\ge 1-\alpha\right\}.
\end{equation}
Note that $ F^P_{Z|\w_1}$ and $\R^P_{|\w_1}(Z)$ are uniquely defined for $\w_{1}\in   \supp (P_1)$, and $\R^P_{|\w_1}(Z)$ can be arbitrary for $\w_{1}\in \O_1\setminus  \supp (P_1)$.

\subsubsection{Multistage setting.}
This can be extended in a straightforward way to the setting where $\O=\O_1\times\cdots\times \O_T$ and $\F=\F_1\otimes \cdots \otimes \F_T$, with  $T\ge 2$. For $P\in \cP$ and $t\in \{2,...,T\}$  we can define conditional probabilities $P_{|\w_{[t-1]}}$ determined by the respective Regular Probability Kernel
$$\bbk_t^P:\F_t\times \O_{[t-1]}\to [0,1],$$
 where $\w_{[t-1]}=(\w_1,...,\w_{t-1})$,
$\O_{[t-1]}=\O_1\times\cdots\times \O_{t-1}$ and $\F_{[t-1]}=\F_1\otimes\cdots\otimes \F_{t-1}$.
Note that
$P_{|\w_{[t-1]}}$ is a probability measure on $(\O_t,\F_t)$ for given (conditional on) $\w_{[t-1]}$.

Then for $Z_t=Z_t(\w_1,...,\w_t)$,
and $$Z_{\w_{[t-1]}}(\cdot):=Z_t(\w_{[t-1]},\cdot)$$
 viewed as a random variable on $(\O_t,\F_t)$,
 define
\begin{equation}\label{condstage-1}
 F^P_{Z_t|\w_{[t-1]}}(z):=P_{|\w_{[t-1]}}(Z_{\w_{[t-1]}}\le z), \;P_t- {\rm  a.e.} \;\;\w_{[t-1]}\in \O_{[t-1]}.
\end{equation}
The corresponding conditional law invariant functional $\R^P_{|\w_{[t-1]}}(Z_t)$ is defined as a function of $F^P_{Z_t|\w_{[t-1]}}$, in the way similar to \eqref{condcdf-2}, that is
\begin{equation}\label{condstage-2}
  \R^P_{|\w_{[t-1]}}(Z_t)=\rho_{t-1}\big(F^P_{Z_t|\w_{[t-1]}}\big).
\end{equation}

 \paragraph{Rectangularity.}
  In particular suppose that $P=P_1\times\cdots\times P_T$. We refer to such setting as {\em  rectangularity} condition.
  In that case $P_{|\w_{[t-1]}}= P_{t}$  for   $P_{[t-1]}$-almost every $\w_{[t-1]}\in \O_{[t-1]}$, where $P_{[t-1]}=P_1\times\cdots\times P_{t-1}$ is the respective marginal distribution on $(\O_{[t-1]},\F_{[t-1]})$.
  Consequently  we can definite the corresponding conditional risk measure as
   \begin{equation}\label{condcdf-4}
\R^P_{|\w_{[t-1]}}(Z_t):=\R^{P_t} (Z_{\w_{[t-1]}}).
\end{equation}
 For example for the Value-at-Risk, definition \eqref{condcdf-3} becomes
 \begin{equation}\label{condcdf-5}
\VaR^P_{\alpha|\w_{[t-1]}}(Z_t)=\inf\left\{z:P_t(Z_{\w_{[t-1]}})\le z)\ge 1-\alpha\right\}.
\end{equation}

Let us consider now robust conditional counterpart of $\bfR$ defined in \eqref{robcoun}.
Suppose that the ambiguity set is  of the following form (the rectangularity assumption)
\begin{equation}\label{rectan-1}
	\cM:=\{P=P_1\times\cdots  P_T: P_t\in \cM_t, \;t=1,...,T\},
\end{equation}
where $\cM_t\subset \cP_t$   are sets of marginal probability   measures.
In that case the robust counterpart of \eqref{condcdf-4} is
   \begin{equation}\label{rectan-2}
\bfR_{|\w_{[t-1]}}(Z_t):=\sup_{P_t\in \cM_t}\R^{P_t}_{|\w_{[t-1]}}(Z_t).
\end{equation}

\begin{remark}
In the rectangular case we can proceed without assuming the conditional law invariance. For example we can use
distributionally robust functionals
\[
\R^{P_t}(\cdot):=\sup_{P_t\in \cM_t}\bbe^{P_t}(\cdot)
\]
even if the ambiguity set $\cM_t$ consists of probability measures which are not absolutely continuous with respect to the reference measure. For example, the ambiguity set $\cM_t$ can consist of probability measures with Wasserstein
distance less than a positive constant  $r$    from a reference probability measure (Wasserstein ball of radius $r$). On the other hand,  in the non-rectangular settings the assumption of conditional law invariance  is essential in the construction of conditional counterparts of risk functionals.
\end{remark}

\paragraph{Nested functionals.} Consider risk functions (risk measures) defined in \eqref{condstage-2}.  For $t=T$ the functional $\R^P_{|\w_{[T-1]}}$ maps random  variable  $Z_T=Z_T(\w_1,...,\w_T)$ into $Z_{T-1}=Z_{T-1}(\w_1,...,\w_{T-1})$, i.e., $Z_{T-1}=\R^P_{|\w_{[T-1]}}(Z_T)$. Recall that $Z_{T}:\O_{[T]}\to\bbr$ is a measurable function and $Z_{T-1}:\O_{[T-1]}\to\bbr$ is measurable by the construction, and that $\O_{[T]}=\O$.
We continue this process iteratively by defining
\begin{equation}\label{iterat-1}
  Z_{t-1}=\R^P_{|\w_{[t-1]}}(Z_t),\;t=T,...,2.
\end{equation}
The corresponding  nested functional is given by $\R^P_{|\w_{1}}(Z_2)$, which is a real number. That is, the nested risk functional is
\begin{equation}\label{iterat-2}
  \cR^P (Z_T)=\R^P_{|\w_{1}} \left (\cdots \R^P_{|\w_{[T-1]}}(Z_T)\right).
\end{equation}
We write this as
  $\cR^P =\R^P_{|\w_{1}}\circ\cdots \circ \R^P_{|\w_{[T-1]}}$.

   In the rectangular case we can use risk functionals of the form \eqref{condcdf-4}.
  For robust functionals of the form \eqref{rectan-2}, the corresponding nested functional ${\bf R}=\bfR_{|\w_{1}}   \circ\cdots \circ \bfR_{|\w_{[T-1]}}:\O\to \bbr$ is defined in the similar way,
\begin{equation}\label{iterat-3}
{\bf R} \left (Z_T\right )=\bfR_{|\w_{1}} \left (\cdots \bfR_{|\w_{[T-1]}}(Z_T)\right).
\end{equation}

   It should be verified that the considered functionals are well defined. As it was already mentioned, $\F_{[t-1]}$-measurability of $\R^P_{|\w_{[t-1]}}(Z_t)$  follows by the construction.    In general, the considered variables $Z_t$ may be restricted to an appropriate linear space  of measurable functions (such linear space is explicitly mentioned in the definition of axioms (A1)-(A4)).
    For the   Value-at-Risk functionals, the corresponding linear space consists of all measurable functions, and $\R^P_{|\w_{[t-1]}}(Z_t)=\VaR^P_{\alpha|\w_{[t-1]}}(Z_t)$ is well defined for any $\F_{[t]}$-measurable $Z_t$.

    For the robust functionals defined in \eqref{rectan-2}, such verification could be more involved. Of course if the space $\O$ is finite, the measurability and existence  issues hold automatically.

\setcounter{equation}{0}
\section{Risk averse Stochastic Optimal Control}
\label{sec-soc}

Consider     Stochastic Optimal Control  (SOC) problem (e.g., \cite{ber78}):
\begin{equation}\label{soc}
\min\limits_{\pi\in \Pi}  \bbe^\pi\left [ \sum_{t=1}^{T}
c_t(x_t,u_t,\xi_t)+c_{T+1}(x_{T+1})
\right],
\end{equation}
where $\xi_1,...,\xi_T$ is a sequence of random vectors, and
$\Pi$ is the set of polices governed by the functional equation
$x_{t+1}=\Phi_t(x_t,u_t,\xi_t),$  $t=1,...,T$, and the  minimization    over controls  $u_t\in \U_t$. In SOC it is usually assumed that   random vectors $\xi_t$  are mutually independent of each other (the stagewise independence assumption), while the distribution of $\xi_t$ is allowed to  depend on state $x_t$ and control $u_t$,  $t=1,...,T$.

Unless stated otherwise, we assume that the distribution of $(\xi_1,...,\xi_T)$
{\em does not}  depend on states and controls, while we allow an interstage dependence of $\xi_1,...,\xi_T$. In that case it suffices to   consider policies
\begin{equation}\label{soc-b}
\Pi=\Big\{\pi=(\pi_1,\ldots,\pi_T):
u_t=\pi_t(x_t,\xi_{[t-1]}),
u_t\in \U_t, x_{t+1}=\Phi_t(x_t,u_t,\xi_t),\;\;t=1,...,T\Big\},
\end{equation}
 where $\xi_{[t]}:=(\xi_1,...,\xi_t)$ denotes history of the process, with (deterministic) initial values  $x_1$ and $\xi_0$.
 We assume that    $\xi_t\in \Xi_t$  where    $\Xi_t\subset \bbr^{d_t}$ is  a closed set, $\F_t$ is the Borel sigma algebra of $\Xi_t$,
 that $\U_t$ is a (nonempty) closed subset of $\bbr^{m_t}$, $\X_t$ is a closed subset of $\bbr^{n_t}$,
 and that the cost functions  $c_t:\X_t\times \U_t\times \Xi_t\to \bbr$ and mappings $\Phi_t:\X_t\times \U_t\times \Xi_t\to \X_{t+1}$ are continuous,   $t=1,...,T$.

It is well known that
    dynamic   equations for value functions of   problem \eqref{soc}  can be written as follows:
 $V_{T+1}(x_{T+1})=c_{T+1}(x_{T+1})$, and for $t=T,...,1$,
 \begin{equation}\label{dyns-1}
V_t(x_t,\xi_{[t-1]})=\inf_{u_t\in \U_t}   \bbe_{|\xi_{[t-1]} }  \big[
c_t(x_t,u_t,\xi_t)+V_{t+1}(\Phi_t(x_t,u_t,\xi_t),\xi_{[t]} ) \big].
\end{equation}
Moreover, if the process $\xi_1,...,\xi_T$  is stagewise independent, then the value functions depend only on state variables, i.e. the dynamic equations take the form
 \begin{equation}\label{dyns-2}
V_t(x_t)=\inf_{u_t\in \U_t}   \bbe  \big[
c_t(x_t,u_t,\xi_t)+V_{t+1}(\Phi_t(x_t,u_t,\xi_t) ) \big].
\end{equation}

Let us consider the following  risk averse counterpart of the SOC.
Consider the set  $\Xi:=\Xi_1\times \cdots\times \Xi_T$ equipped with  sigma algebra $\F=\F_1\otimes \cdots \otimes \F_T$, and denote by $\cP$ the set of probability measures on $(\Xi,\F)$. For $P\in \cP$ consider conditional law invariant risk   functionals $\R^P_{|\xi_{[t-1]}}$, defined as in
\eqref{condstage-2}, and the corresponding nested functional $\cR^P =\R^P_{|\xi_{1}}\circ\cdots \circ \R^P_{|\xi_{[T-1]}}$, defined as in \eqref{iterat-2}. The respective risk averse counterpart  of problem \eqref{soc} is
\begin{equation}\label{socrisk-1}
\min\limits_{\pi\in \Pi}  \cR^P \left [ \sum_{t=1}^{T}
c_t(x_t,u_t,\xi_t)+c_{T+1}(x_{T+1})
\right],
\end{equation}
where the minimization is over policies of the form \eqref{soc-b}.
The dynamic programming equations for problem \eqref{socrisk-1} can be written as
$V_{T+1}(x_{T+1})=c_{T+1}(x_{T+1})$, and for $t=T,...,1$,
 \begin{equation}\label{socrisk-2}
V_t(x_t,\xi_{[t-1]})=\inf_{u_t\in \U_t}   \R^P_{|\xi_{[t-1]}}  \big[
c_t(x_t,u_t,\xi_t)+V_{t+1}(\Phi_t(x_t,u_t,\xi_t),\xi_{[t]} ) \big],
\end{equation}
which  can be viewed as the counterpart of   dynamic equations \eqref{dyns-1}.
A sufficient condition for policy $\bar{u}_t=\bar{\pi}_t(x_t,\xi_{[t-1]})$, $t=1,...,T$, to be optimal is that
 \begin{equation}\label{socrisk-3}
\bar{u}_t\in \argmin_{u_t\in \U_t}   \R^P_{|\xi_{[t-1]}}  \big[
c_t(x_t,u_t,\xi_t)+V_{t+1}(\Phi_t(x_t,u_t,\xi_t),\xi_{[t]} ) \big].
\end{equation}
Proof of the above dynamic equations is based on an  interchangeability property  of risk functionals  satisfying the axioms of monotonicity and translation equivariance (cf., \cite{shaope2017}). It is also pointed out in  \cite{shaope2017} that without the {\em strict} monotonicity property of the risk functionals, conditions \eqref{socrisk-3} could be not necessary for  optimality, i.e. it could exist a  policy $\bar{\pi}\in \Pi$  which is optimal for problem \eqref{socrisk-1} but does not satisfy dynamic equations \eqref{socrisk-3}. Such optimal  policy is not dynamically consistent in the sense that although it is optimal from the point of view of the first stage, it could be not optimal at the later stages conditional  on  some realizations of the data process.

\paragraph{Rectangular setting.}
Consider now the rectangular setting, i.e., assume that $P=P_1\times \cdots\times P_T$, where $P_t$ is a marginal distribution of $\xi_t$.
  Then risk functionals become of the form   \eqref{condcdf-4}, that is
$\R^P_{|\xi_{[t-1]}}(Z_t):=\R^{P_t} (Z_{\xi_{[t-1]}})$. In that case value functions depend on state variables only and dynamic equations become
  \begin{equation}\label{dynrobust}
V_t(x_t)=\inf_{u_t\in \U_t}   \R^{P_t}   \big[
c_t(x_t,u_t,\xi_t)+V_{t+1}(\Phi_t(x_t,u_t,\xi_t) ) \big],
\end{equation}
where $\xi_t \sim  P_t$, $t=1,...,T$.
The rectangularity condition is a counterpart of stagewise independence condition,  and risk averse dynamic equations \eqref{dynrobust} is the  counterpart  of equations \eqref{dyns-2}.

    In the rectangular case we can also consider the robust risk averse setting. That is, let $\cM_t$ be a set of probability measures on $(\Xi_t,\F_t)$. Then counterparts of the dynamic equations \eqref{dynrobust} become
     \begin{equation}\label{dynrobust-2}
V_t(x_t)=\inf_{u_t\in \U_t} \sup_{P_t\in \cM_t}  \R^{P_t}   \big[
c_t(x_t,u_t,\xi_t)+V_{t+1}(\Phi_t(x_t,u_t,\xi_t) \big].
\end{equation}
     This corresponds to replacing functional $\R^{P_t}(\cdot)$ with its robust counterpart
      (compare with \eqref{robcoun})
      \begin{equation}\label{robcoun-d}
     \bfR_t(\cdot):=\sup_{P_t\in \cM_t} \R^{P_t}(\cdot),\;t=1,...,T.
      \end{equation}
These dynamic equations correspond to the nested formulation of the risk averse SOC problem:
 \begin{equation}\label{socriskrobust}
\min\limits_{\pi\in \Pi}  {\bf R} \left [ \sum_{t=1}^{T}
c_t(x_t,u_t,\xi_t)+c_{T+1}(x_{T+1})
\right],
\end{equation}
 where
  ${\bf R}=\bfR_{|\xi_{1}}   \circ\cdots \circ \bfR_{|\xi_{[T-1]}}$ (compare with \eqref{iterat-3}).

\subsection{Infinite horizon risk averse SOC.}
\label{sec-infsoc}
  We also can consider risk averse infinite horizon SOC problem  and its robust  extension.
 The risk averse  counterpart of the classical    Bellman equation is
  \begin{equation}\label{bellman_soc_risk_averse}
V(x)=\inf_{u \in \U }   \R^{P}   \big[
c (x,u,\xi)+\beta V(\Phi(x,u,\xi) ) \big],
 \end{equation}
 and
   \begin{equation}\label{bellrobust}
V(x)=\inf_{u \in \U } \sup_{P\in \cM}  \R^{P}   \big[
c (x,u,\xi)+\beta V(\Phi(x,u,\xi) ) \big],
 \end{equation}
 for the robust risk averse problem.
 Here $\beta\in (0,1)$  is the discount factor, $\X\subset \bbr^n$, $\U\subset \bbr^m$,  $\Xi\subset \bbr^d$    are  (nonempty) closed sets, $\F$ is the Borel sigma algebra of $\Xi$,  $c:\X\times \U\times \Xi\to \bbr$ and $\Phi:\X\times \U\times \Xi\to \X$ are continuous functions, and
  $\cM$ is a set of probability distributions on $(\Xi,\F)$.

\begin{itemize}
  \item []
 We assume that for every $P\in \cM$,  the functional $\R^P$ satisfies the axioms of monotonicity and translation equivariance, and that  $\R^P(0)=0.$
\end{itemize}

It follows that   if $|Z|\le \kappa$,  $P$-almost surely for some $\kappa\in \bbr$, then $|\R^P(Z)|\le    \R(\kappa)=\R(0)+\kappa=\kappa$.
 Assuming that the cost function is bounded, it follows   that the corresponding Bellman operator has the contraction property, and hence by the fixed point theorem,  equation \eqref{bellrobust}  has a unique solution.
 That is, consider the space $\bbb$ of bounded measurable functions $g:\X\to \bbr$, equipped with the sup-norm $\|g\|_\infty=\sup_{x\in \X}|g(x)|$, and Bellman operator $\T:\bbb\to\bbb$ ,
\begin{align}
\label{belloper}
\T(g)(\cdot):=
\inf_{u \in \U}  \sup_{P\in \cM} \R^{P}   \big[
c(\cdot ,u ,\xi )+ \beta g(\Phi (\cdot,u ,\xi) ) \big],\;g\in \bbb.
\end{align}
We have   that  if $|Z|\le \kappa$,  $P$-almost surely for some $\kappa\in \bbr$, then $|\R^P(Z)|\le  \kappa$. This implies that  if  the cost function is bounded,  then  for any  bounded  $g:\X\to\bbr$, it follows that  $\T(g)$ is bounded. That is, the operator $\T$  maps the space of bounded functions into   the space of bounded functions.

  The operator $\T$ has the following properties. It is monotone, i.e., if $g,g'\in \bbb$ and $g(\cdot)\ge g'(\cdot)$, then $\T(g)(\cdot)\ge \T(g')(\cdot)$. This follows from the monotonicity property of $\R^P$. Also   from the translation equivariance property of  $\R^P$ follows   the  property of constant shift: for any $g\in \bbb$ and $c\in \bbr$,
  $\T(g+c)=\T(g)+\beta c$. It follows that $\T$ is a contraction mapping, i.e.,
  \[
  \|\T(g)-\T(g')\|_\infty \le \beta \|g-g'\|_\infty,\;g,g'\in \bbb,
\]
and hence by the Banach Fixed Point Theorem has unique fixed point, denoted by $V$.
Equation \eqref{bellrobust} corresponds to the nested formulation of the respective infinite horizon risk averse problem.

\begin{remark}
To ensure that \eqref{bellrobust} is well-defined in the first place, we proceed as follows.
Suppose that the set $\Xi$ is compact and consider the space $C(\Xi)$ of continuous functions $\phi:\Xi\to\bbr$ equipped with the sup-norm. The dual  $C(\Xi)^*$ of that space is the space of finite signed (Borel) measures on $\Xi$ (Riesz representation).
Since $\Xi$ is compact,   the weak$^*$ topology  of $C(\Xi)^*$  is metrizable.
Let $\cM$ be weakly$^*$-closed (and hence weakly$^*$-compact by the Banach - Alaoglu theorem) and $\U$ be compact.
Suppose $\R^P(Z)$ is continuous in $P$ in weak$^*$-topology, and continuous in $Z$ in the sup-norm, then
it is clear that for any continuous $g(\cdot)$, the optimization problem on the right hand side of \eqref{belloper} is continuous in $(u, x, P)$.
Since $\cM$ and $\U$ are compact, it follows that $\T(g)(\cdot)$ is continuous.
Combining this with $\T$ being contractive, it follows that the value function $V$ is continuous and \eqref{bellrobust} is well-defined.
\end{remark}

\subsubsection{Sample complexity of infinite horizon risk averse SOC}

We focus in this section on sample complexity in  the infinite-horizon setting,    the sample complexity for the finite-horizon setting can be derived similarly. Specifically we deal with the Value-at-Risk measure, and the analysis    can be viewed as an extension of   section
\ref{sec_sample_var}.

 Let $\xi_1,...,\xi_N$  be iid samples of $\xi\sim P$,  and $P_N=N^{-1}\sum_{i=1}^N \delta_{\xi_i}$ be the corresponding empirical measure. The   empirical counterpart of the Bellman equation is obtained by replacing probability measure $P$ in \eqref{bellman_soc_risk_averse} with the empirical measure $P_N$.
Consider the counterpart of Bellman operator \eqref{belloper} with respect to $P_N$, that is
\begin{align*}
\T_N(g)(\cdot):=
\inf_{u \in \U}   \R^{P_N}   \big[
c(\cdot ,u ,\xi )+g(\Phi (\cdot,u ,\xi) ) \big].
\end{align*}
Assuming  that  the cost function is bounded, the operator $T_N$ has a unique fixed point denoted  $V_N$, which is the solution of the empirical counterpart of the  Bellman equation. Our interest in this section is to determine the sample size $N$ such that $  \| V - V_N \|_\infty<\e$ with   high probability  for a given $\e>0$.
We make the following analogue of Assumption \ref{ass-lip}.

\begin{assumption}
\label{ass-lip_soc-1}
(i) The sets $\X$ and $\U$ are compact. (ii)
There is  a positive constant $L$  such that
\begin{eqnarray*}
&&\left| c(x, u,\xi) - c(x', u',\xi) \right|    \leq L \| (x, u) - (x', u') \|,\\
&&  \left\| \Phi(x, u,\xi) - \Phi(x', u',\xi) \right\|    \leq L \| (x, u) - (x', u') \|,
\end{eqnarray*}
for all $x, x' \in \X$, $u, u' \in    \U$,  $\xi \in   \Xi$.
\end{assumption}
Unfortunately the above assumption does not guarantee Lipschitz continuity of the value function $V(\cdot)$. Nevertheless we can proceed as follows.
Let $\tilde{V}$ denote an approximation of $V$ such that $\tilde{V}$ is $\tilde{L}$-Lipschitz continuous.
Then
\begin{align*}
\| \tilde{V} - V_N \|_\infty  & = \| \T \tilde{V} + \tilde{V} - \T \tilde{V} - \T_N V_N \|_\infty  \\
& \leq \| \T \tilde{V} - \T_N V_N \|_\infty + \| \tilde{V} - \T \tilde{V} \|_\infty \\
& =   \| \T \tilde{V} - \T_N \tilde{V} + \T_N \tilde{V} -  \T_N V_N \|_\infty + \| \tilde{V} - \T \tilde{V} \|_\infty \\
& \leq \| (\T - \T_N) \tilde{V} \|_\infty +  \beta \| \tilde{V} - V_N \|_\infty + \| \tilde{V} - \T \tilde{V} \|_\infty ,
\end{align*}
and consequently
\begin{align}
 \label{soc_saa_gap_conversion}
\| \tilde{V} - V_N \|_\infty
\leq \tfrac{1}{1-\beta} \big(
\| (\T - \T_N) \tilde{V} \|_\infty
+ \| \tilde{V} - \T \tilde{V} \|_\infty
\big).
\end{align}
This immediately implies
\begin{align}
\| V - V_N  \|_\infty
& \leq
\| V - \tilde{V} \|_\infty + \| \tilde{V} - V_N \|_\infty \nonumber  \\
& \leq
\| V - \tilde{V} \|_\infty  +
\tfrac{1}{1-\beta} \big(
\| (\T - \T_N) \tilde{V} \|_\infty
+ \| \tilde{V} - V + \T V -  \T \tilde{V} \|_\infty
\big)  \nonumber \\
& \leq \tfrac{1}{1-\beta} \big( \| (\T - \T_N) \tilde{V} \|_\infty
+ 2 \| V - \tilde{V} \|_\infty \big). \label{diff_v_v_N_soc}
\end{align}


\begin{lemma}\label{lemma_lipschitz_approx_v}
Suppose that Assumption \ref{ass-lip_soc-1} holds and $\R^P(Z)$ is $\LL_{\R}$-Lipschitz continuous in $Z$ w.r.t. $\|\cdot \|_\infty$ norm,
i.e.,
\begin{align}\label{condition_lipschitz_risk}
|\R^P(Z) - \R^P(Z')| \leq \LL_{\R} \| Z - Z' \|_\infty.
\end{align}
Then for any $\e > 0$, there exists $\tilde{V}$ such that
$\| \tilde{V} - V \|_\infty \leq \e$, and
$\tilde{V}$ is $\tilde{L}$-Lipschitz continuous with
\begin{equation}\label{lipschitz_constant_approx_v}
\tilde{L} := \tfrac{ \LL_{\R} L [ (\beta \LL_{\R} L)^k - 1 ] }{\beta \LL_{\R} L - 1} \; \;and\;\;
k := \tfrac{1}{1-\beta} \log \big(\tfrac{1}{\e} \big) ,
\end{equation}
where  we   assume (without loss of generality) that $\beta \LL_{\R} L > 1$.
\end{lemma}

\begin{proof}
Consider a  function $g:\X\to \bbr $ that is $L_g$-Lipschitz continuous, and define
\begin{align*}
Q_g (x, u, \xi) = c(x, u, \xi) + \beta g( \Phi(x, u, \xi)).
\end{align*}
Then we have $\T(g)(x) = \inf_{u \in \U} \R^P(Q_g(x, u, \cdot))$.  Furthermore
\begin{eqnarray*}
 \left| \R^P(Q_g(x, u, \cdot)) -  \R^P(Q_g(x', u', \cdot)) \right|
&\leq&
\LL_{\R} \| Q_g(x, u, \cdot)  -  Q_g(x', u', \cdot) \|_\infty\\
&\leq &  \LL_{\R} L  \left(
1 + \beta L_g
\right)
\| (x, u) - (x', u') \|,
\end{eqnarray*}
 from which we conclude that
$\T(g)(\cdot)$ is also $\LL_{\R} L  \left(
1 + \beta L_g
\right)$-Lipschitz continuous.

Now consider $V^{(0)} \equiv 0$, and define $V^{(k)} = \T (V^{(k-1)})$.
Denote $L_k$ the Lipschitz constant of $V^{(k)}$, then it is clear that $L_k= 0$ and
\begin{align*}
L_{k} \leq \LL_{\R} L + \beta \LL_{\R} L L_{k-1}
\leq
\begin{cases} \tfrac{ \LL_{\R} L [ (\beta \LL_{\R} L)^k - 1 ] }{\beta \LL_{\R} L - 1}, & \beta \LL_\R L \neq 1, \\
k \LL_\R L , & \beta \LL_\R L = 1,
\end{cases}
\end{align*}
for any $k \geq 1$.
Since $\T(\cdot)$ is a contractive operator in $\| \cdot \|_\infty$-norm, it suffices to take $k= \tfrac{1}{1-\beta} \log(\tfrac{1}{\e})$ to ensure
$\| V - V^{(k)} \|_\infty \leq \e$.
Taking $\tilde{V} = V^{(k)}$ concludes the proof.
\end{proof}

In the remainder of this section we deal with
$\R^P(\cdot) := \VaR^P_\alpha (\cdot)$.
From the first inequality in \eqref{condition_lipscthiz_in_x_sample-1}, it is clear that in the case of  $\R^P(\cdot) := \VaR^P_\alpha (\cdot)$   we  can take $\LL_{\R} = 1$ in \eqref{condition_lipschitz_risk}.

Consider now the case where the probability measure  $P$ has finite support. For  $\alpha \in (0,1)$, define
\begin{eqnarray}
\nonumber
& l_\alpha := \sup \{    P(A): P(A) \leq 1 - \alpha,\;A\in \F \},\; \;r_\alpha := \inf \{   P(A): P(A) \geq 1 - \alpha,\;A\in \F \},\\
&  \kappa_\alpha := \min \{ (1- \alpha) - l_\alpha, r_\alpha - (1- \alpha ) \}.
\label{leftright}
\end{eqnarray}
{Since $P$ has finite support}, the set of $\alpha$ where $\kappa_\alpha = 0$ is finite.
The condition $\kappa_\alpha>0$ and the following theorem  can be viewed as a counterpart of condition \eqref{epsilon} and Proposition \ref{prop_cvar_saa_exact}. The estimate \eqref{estsamp}, below,   suggests that in the present case  the
number of samples needed in the risk averse setting only grows  linearly (up to an logarithmic factor) in $(1-\beta)^{-1}$ as the  discount factor $\beta$ approaches $1$.

\begin{theorem}\label{thrm_sample_soc_discrete}
Let $\R^P = \VaR_\alpha^P$  in \eqref{bellman_soc_risk_averse}, for some   $\alpha \in (0,1)$.
Suppose Assumption {\rm \ref{ass-lip_soc-1}} holds and the probability measure  $P$ has finite support. Suppose further that  $\kappa_\alpha  > 0$,  with  $\kappa_\alpha$ defined in \eqref{leftright}.
  Let  $D$ denote the diameter of the set $\X \times \U$. Then
for any $\delta \in (0,1)$,  $\e > 0$, and
\begin{align}
\label{estsamp}
N \geq
   \tfrac{1}{2}  \kappa_\alpha^{-2} \left[ (n+m) \log ( \tfrac{ 8  D  L^2}{ \e (1-\beta) (\beta  L -1)  } ) + \tfrac{n+m }{1-\beta} \log(\beta  L) \log (\tfrac{4}{\e (1-\beta)})  + \log(\tfrac{2}{\delta}) \right],
\end{align}
it follows that $\| V - V_N \|_\infty \leq \e$
with probability at least $1-\delta$.
\end{theorem}

\begin{proof}
To begin, let us recall Lemma \ref{lemma_lipschitz_approx_v}, and let $\tilde{V}$ denote the approximation of $V$ such that
\begin{align}\label{approx_v_by_tilde_v}
\| V - \tilde{V} \|_\infty \leq (1-\beta) \e /4
\end{align}
with Lipschitz constant
$
\tilde{L} = \tfrac{  L [ (\beta  L)^k - 1 ] }{\beta  L - 1}
$,
$ k= \tfrac{1}{1-\beta} \log \big(\tfrac{4}{\e (1-\beta)} \big).
$
 For notational convenience, define
\begin{align*}
\tilde{Z}_{x, u} (\xi) := c(x ,u ,\xi )+ \beta \tilde{V}(\Phi (x ,u ,\xi) )
\end{align*}
and
$
\overline{\T}(\tilde{V})(x, u):=
   \VaR^{P}_\alpha   \big[
\tilde{Z}_{x, u}  \big]
$,
and
$
\overline{\T}_N (\tilde{V})(x, u):=
   \VaR^{P_N}_\alpha   \big[
\tilde{Z}_{x, u} \big].
$
Note that from Assumption \ref{ass-lip_soc-1} and $\tilde{L}$-Lipschitz continuity of $\tilde{V}$, we have
\begin{align}\label{lipschitz_of_Z}
\| \tilde{Z}_{x, u} - \tilde{Z}_{x' , u'} \|_\infty \leq  \overline{L} \| (x,u) - (x', u') \| , ~ \overline{L} =
L (1 + \beta \tilde{L}) .
\end{align}

For a fixed $x, u \in \X \times \U$, from Proposition \ref{prop_cvar_saa_exact}  and the definition of $\kappa_\alpha$,
for $N\ge \half \kappa_\alpha^{-2}\log(2/\delta)$,  with probability $1-\delta$ we have
\begin{align*}
\overline{\T}(\tilde{V})(x, u) = \overline{\T}_N (\tilde{V})(x, u).
\end{align*}
Consider an $\eta$-net of ${\cal X} \times \U$, denoted by ${\cN}_{\eta}$. That is,
for any $(x, u) \in {\cal X} \times \U$, there is  $(x_\eta, u_\eta) \in {\cN}_{\eta}$ such that $\left\Vert (x, u) - (x_\eta, u_\eta) \right\Vert \leq \eta$.
For $(x, u) \in \X \times \U$ we can write
\begin{equation}
\label{ineq-2}
\begin{array}{lr}
&\left| \VaR_\alpha^{P} \left (\tilde{Z}_{x, u} \right) - \VaR_\alpha^{P_N} \left (\tilde{Z}_{x, u} \right) \right|
 \leq
\left| \VaR_\alpha^{P} \left (\tilde{Z}_{x_\eta, u_\eta}) \right) - \VaR_\alpha^{P_N} \left (\tilde{Z}_{x_\eta, u_\eta}  \right) \right| +  \\
 & \left |\VaR_\alpha^{P} \left (\tilde{Z}_{x, u} \right) - \VaR_\alpha^{P} \left (\tilde{Z}_{x_\eta, u_\eta} \right) \right|
  +
\left| \VaR_\alpha^{P_N} \left (\tilde{Z}_{x, u}   \right) - \VaR_\alpha^{P_N} \left ( \tilde{Z}_{x_\eta, u_\eta}  \right) \right|.
\end{array}
\end{equation}
Hence for $N\ge \half \kappa_\alpha^{-2}\log(2 |\cN_\eta| /\delta)$ and $\eta \leq \e (1-\beta) / (4\overline{L})$, with probability at least $1-\delta$, for any $(x, u) \in \X \times \U$ we have
\begin{align*}
& \left| \overline{\T}(\tilde{V})(x, u)  - \overline{\T}_N(\tilde{V})(x, u)  \right| \\
\leq &
\left| \VaR_\alpha^{P} \left (\tilde{Z}_{x, u} \right) - \VaR_\alpha^{P} \left (\tilde{Z}_{x_\eta , u_\eta} \right) \right|
+
\left| \VaR_\alpha^{P_N} \left (\tilde{Z}_{x, u} \right) - \VaR_\alpha^{P_N} \left (\tilde{Z}_{x_\eta , u_\eta} \right) \right| \\
  \overset{(a)}{\leq} &
2 \overline{L} \eta \leq (1-\beta) \e/2,
\end{align*}
where
$(a)$ follows from $\VaR_\alpha^P(\cdot)$ being $1$-Lipschitz continuous w.r.t $\| \cdot \|_\infty$-norm and \eqref{lipschitz_of_Z}.
Since $\T (\tilde{V}) (x)  = \inf_{u \in \U} \overline{\T}(\tilde{V}) (x, u)$
and $\T_N (\tilde{V}) (x)  = \inf_{u \in \U} \overline{\T}_N(\tilde{V}) (x, u)$,
 this  immediately implies
$$
\| (\T - \T_N) \tilde{V} \|_\infty \leq (1-\beta) \e/2.
$$
Combining this together with \eqref{approx_v_by_tilde_v} and   \eqref{diff_v_v_N_soc}  yields $\| V -  V_N \|_\infty\leq \e$.
The rest of the proof pertains to bounding $| \cN_\eta|$ and follows similar argument as in Theorem \ref{pr-unifbound}.
\end{proof}

 \begin{remark}
For continuous $P$, following similar arguments as in Theorem \ref{pr-unifbound} and \ref{thrm_sample_soc_discrete}, a sample complexity of ${\cal O} ((1-\beta)^{-2} c^{-2} \e^{-2}) $ could be established if $\beta   \in (0, L^{-1})$, and there are constants $b>0$ and ${c} >0$ such that for all $(x,u)\in \X\times \U$,
 \begin{equation}\label{growth-2a_soc}
  F_{x, u} (z')-F_{x, u}(z)\ge  {c} \, (z'-z) \;\text{for all}\; z,z' \in [\nu_{x, u} -b,\nu_{x, u} +b]\;\text{such that}\; \;z'\ge z,
\end{equation}
 where $F_{x, u}$ is the cdf of $Z_{x, u} (\xi) : =  c(x ,u ,\xi )+ {V}(\Phi (x ,u ,\xi) )$ and
 $\nu_{x, u}:=\VaR^{P}_{\alpha}(Z_{x, u})$.
That is,  the sample size, required to attain a given
{\em relative} error of the empirical (SAA) solution, is not sensitive to the discount factor, even if the discount factor is very
close to one.
 This is in line with the similar observation made in \cite{cltshapcheng21}, where it was based on different arguments.
It should be noted, however,  that the above  \eqref{growth-2a_soc} depends on the unknown value function $V(\cdot)$ and is in general difficult to verify.
In particular, it remains open to establish the existence of density function\footnote{
When $\beta L < 1$, one can show that $V$ is Lipschitz continuous, and hence $Z_{x,u}(\xi)$ is Lipschitz continuous in $\xi$.
Suppose $P$ is absolutely continuous w.r.t. the Lebesgue measure, and the first-order stationary points of $Z_{x,u}(\cdot)$ take Lebesgue measure zero, then from   Coarea formula  \cite{Federer},  it follows that $Z_{x,u}^{-1}(A)$ is a $P$-null set for any $A \subset \mathbb{R}$ taking Lebesgue measure zero.
In this case, the density function of $Z_{x,u}$ exists due to Radon-Nikodym theorem.} for $Z_{x,u}$.
 \end{remark}

\subsection{Randomized policies}
Consider the (robust)  rectangular setting and suppose that the decision maker can choose controls at random. That is, for $t=1,...,T,$  let   $\cS_t$ be  the set of (Borel) probability measures on $\U_t$,  and let $(u_t,\xi_t)\sim Q_t\times P_t$  be random with $Q_t\in \cS_t$ and $P_t\in \cM_t$.
Consider the following extension of dynamic equations \eqref{dynrobust-2}: $\sV_{T+1}(x_{T+1})=c_{T+1}(x_{T+1})$, and for $t=T,...,1$,
    \begin{equation}\label{dynrandom-1}
\sV_t(x_t)=\inf_{Q_t\in \cS_t}   \sup_{P_t\in \cM_t}   \R^{Q_t\times P_t}   \big [
c_t(x_t,u_t,\xi_t)+\sV_{t+1}(\Phi_t(x_t,u_t,\xi_t) \big ],
\end{equation}
where $\R^{Q_t\times P_t} $ is a risk functional on a linear space of random variables $Z:\U_t\times \Xi_t\to \bbr$. Specifically we assume that
   \begin{equation}\label{dynrandom-2}
 \R^{Q_t\times P_t} (Z_t(u_t,\xi_t)):=\bbe^{Q_t}\left[\R^{P_t}(Z_t(u_t,\xi_t))\right],
\end{equation}
where
$u_t\sim Q_t$ and $\R^{P_t}$ is a risk functional with respect to $\xi_t\sim P_t$. In particular if $\R^{P_t}=\bbe^{P_t}$, equation \eqref{dynrandom-1}  becomes
    \begin{equation}\label{dynrandom-3}
\sV_t(x_t)=\inf_{Q_t\in \cS_t}   \sup_{P_t\in \cM_t}    \bbe^{Q_t\times P_t}   \big [ c_t(x_t,u_t,\xi_t)+\sV_{t+1}(\Phi_t(x_t,u_t,\xi_t) \big].
\end{equation}

  We say that there exists a non-randomized optimal policy if for $t=1,...,T$,  the minimum in the right hand side of  \eqref{dynrandom-1} is  attained at   Dirac measure $Q_t=\delta_{u_t}$. That is, existence of the non-randomized policy means that dynamic equations \eqref{dynrandom-1} are equivalent to dynamic equations \eqref{dynrobust-2}.

  The dual the min-max problem in the right hand side of \eqref{dynrandom-1} is obtained by interchanging the `inf' and `sup' operators. That is,  value functions of the dual problem are defined as $\sW_{T+1}(x_{T+1})=c_{T+1}(x_{T+1})$, and for $t=T,...,1$,
     \begin{equation}\label{dyndual-1}
\sW_t(x_t)= \sup_{P_t\in \cM_t} \inf_{Q_t\in \cS_t}  \underbrace{ \bbe^{Q_t} \big [ \R^{P_t} }_{\R^{Q_t\times P_t}}\big ( c_t(x_t,u_t,\xi_t)+\sW_{t+1}(\Phi_t(x_t,u_t,\xi_t) \big)\big ].
\end{equation}
 The function
 \[
 \psi(x_t,u_t):=\R^{P_t}\big ( c_t(x_t,u_t,\xi_t)+\sW_{t+1}(\Phi_t(x_t,u_t,\xi_t) \big)
 \]
  inside the brackets in the right hand side of  \eqref{dyndual-1} is a function of $x_t$ and $u_t$, and consequently    the expectation  $\bbe^{Q_t}[ \psi(x_t,u_t)]$
is minimized over  $Q_t\in \cS_t$.
Let us observe that  for every $x_t$ it suffices  to perform
this minimization
 over Dirac measures $Q_t=\delta_{u_t}$,    $u_t\in  \U_t$. That is, $\sW_t(\cdot)=W_t(\cdot)$, where $W_t$ is  defined by  equation
     \begin{equation}\label{dyndual-2}
W_t(x_t)= \sup_{P_t\in \cM_t} \inf_{u_t\in \U_t}       \R^{P_t} \big ( c_t(x_t,u_t,\xi_t)+W_{t+1}(\Phi_t(x_t,u_t,\xi_t) \big).
\end{equation}
By the standard theory of min-max problems we have that $V_t(\cdot)\ge W_t(\cdot)$. It is said that there is no duality  gap between the    dual problems \eqref{dynrobust-2} and \eqref{dyndual-2} if $V_t(\cdot)=W_t(\cdot)$.

  For $x_t\in \X_t$, it is said that $(u_t^*,P^*_t)$ is a saddle point of  the min-max problem \eqref{dynrobust-2}   if
  \begin{eqnarray}
  \label{saddle-1}
    u^*_t \in  \argmin_{u_t\in \U_t}      \R^{ P^*_t}   \big [
c_t(x_t,u_t,\xi_t)+V_{t+1}(\Phi_t(x_t,u_t,\xi_t) \big ],\\
  \label{saddle-2}
      P^*_t \in \argmax_{P_t\in \cM_t}      \R^{ P_t}   \big [
c_t(x_t,u^*_t,\xi_t)+V_{t+1}(\Phi_t(x_t,u^*_t,\xi_t) \big ].
  \end{eqnarray}
Existence of the saddle point is a sufficient condition for the no duality gap of the   min-max problems   \eqref{dynrobust-2} and \eqref{dyndual-2}. Conversely,  given optimal solutions $u^*_t$ and $P^*_t$ of the respective problems \eqref{dynrobust-2} and \eqref{dyndual-2}, the no duality gap property implies that
$(u^*_t,P^*_t)$ is a saddle point.

In order to verify the no duality gap between the min-max problems in \eqref{dynrandom-1} and \eqref{dyndual-1}, i.e., that $\sV_t(\cdot)=\sW_t(\cdot)$,  we need the following condition.

\begin{assumption}
\label{ass-sion}
For $t=1,...,T,$  the set $\cM_t$ is convex, the set $\U_t$ is compact, and for every  $P\in \cM_t$ the functional $\R^{P}$ is  concave in $P$,  i.e.,  if $P, P'\in \cM_t$ and $\tau\in [0,1]$, then
$
\R^{\tau P+(1-\tau)P'}(\cdot )\ge \tau \R^P(\cdot)+(1-\tau)\R^{P'}(\cdot).
$
\end{assumption}

Suppose that  the above assumption holds. Then  the functional
 $\R^{Q_t\times P_t}$,  defined in \eqref{dynrandom-2}, is concave in $P_t$, and 
 is convex (linear) in $Q_t$.
  By Sion's theorem \cite{sion}, there is no duality gap between  the min-max problems
  \eqref{dynrandom-1} and \eqref{dyndual-1}.
 There is a large class of   risk functionals which are concave with respect to the probability distribution. For example consider   risk functionals of the from $\R^P(Z)=\inf_{\theta\in \Theta}\bbe^P[\Psi(Z,\theta)]$, where $\Theta$  is a subset of a finite dimensional vector space and
$\Psi:\bbr\times \Theta\to \bbr$ is a real valued function. The $\avr_\alpha$ risk measure is of that form.  Since $\R^P(Z)$ is given by the infimum of  linear in $P$ functionals, it is concave in $P$.
 Another example of concave in $P$ functional   is   the Value-at-Risk measure $\R^P=\VaR^P_{\alpha}$.

The following result about   existence of the non-randomized optimal policies is an extension of \cite[Theorem 4.1]{ShapiroLi2025}.

\begin{proposition}
 \label{pr-nonraand}
 Suppose that
for  $t=1,...,T$, there exists saddle point $(u^*_t,P^*_t)$ satisfying conditions \eqref{saddle-1} and \eqref{saddle-2}. Then there exists  non-randomized optimal  policy given by   $\pi^*_t(x_t)=u^*_t$.

Conversely,  suppose that Assumption \ref{ass-sion} holds and there exists non-randomized optimal policy. Then the saddle point exists.
\end{proposition}

\begin{proof}
Suppose that there exists saddle point $(u^*_t,P^*_t)$. Then $V_t(\cdot)=W_t(\cdot)$. Moreover,  we have that   $V_t(\cdot)\ge \sV_t(\cdot)$, and by duality $\sV_t(\cdot)\ge \sW_t(\cdot)$.   Since  $W_t(\cdot)= \sW_t(\cdot)$, it follows that $V_t(\cdot)= \sV_t(\cdot)$. That is, the value functions with respect to randomized and non-randomized policies are the same. This implies existence of optimal policy, which is given by $\pi^*_t(x_t)=u^*_t$.

 Conversely,  suppose  there exists the non-randomized optimal policy. This implies that $V_t(\cdot)= \sV_t(\cdot)$. Moreover, suppose  that Assumption \ref{ass-sion} holds. Then
  by Sion's theorem $\sV_t(\cdot)=\sW_t(\cdot)$.  Since $\sW_t(\cdot)=W_t(\cdot)$, it follows that $V_t(\cdot)= W_t(\cdot)$. Because both respective problems have optimal solutions, existence of the saddle point  follows.
\end{proof}

In particular, suppose that $\cM_t=\{P^*_t\}$ is the singleton and the minimization   problem in  \eqref{saddle-1} has an optimal solution $u_t^*$.
Then clearly  the saddle point exists   and consequently there exists the  non-randomized optimal policy.

Suppose that for $t=1,...,T$,  functionals $\R^{P_t}$ are convex, i.e., satisfy the axiom of convexity. Suppose further that the sets $\U_t$ and  $\cM_t$ are convex,  the cost functions $c_t(x_t,u_t,\xi_t)$ are convex in $(x_t,u_t)$  and the mappings $\Phi_t(x_t,u_t,\xi_t)=A_t(\xi_t) x_t +B_t(\xi_t) u_t+b_t(\xi_t)$ are affine. Then the value functions $V_t(x_t)$ are convex and the min-max problem  in the right hand side of \eqref{dynrobust-2} becomes  convex-concave. Suppose further that the sets $\U_t$ are  compact.
Then   by Sion's theorem there is no duality gap between problems \eqref{dynrobust-2} and \eqref{dyndual-2}, and thus the saddle point exists provided the primal \eqref{dynrobust-2} and dual  \eqref{dyndual-2}  problems have optimal solutions. Consequently in that case existence of non-randomised policies follows (cf., \cite[Section 4]{ShapiroLi2025}).
For non-convex functionals $\R^{P_t}$,  verification of existence of the saddle points could be more involved.

\setcounter{equation}{0}
\section{Risk averse Markov Decision Processes}
\label{sec-mdp}
 Consider a (finite horizon) Markov Decision Process (MDP) (e.g., \cite{puterman2014markov})
\begin{equation}\label{mdprob-1}
 \min_{\pi\in \Pi} \bbe^{\pi}
 \left [ \tsum_{t=1}^{T}
c_t(s_t,a_t,s_{t+1})+c_{T+1}(s_{T+1})
\right].
\end{equation}
Here $\s_t$ is the state space,  $\A_t$ is the action set and $c_t:\s_t\times \A_t\times \s_{t+1}\to \bbr$ is the cost function,
at stage $t=1,...,T$. We assume that the state   $\s_t$ and action $\A_t$ spaces are   Polish spaces, and denote by $\F_t$ the (Borel) sigma algebra of $\s_t$.  The dynamic process is defined by transitional kernels (transitional probabilities)  $P_t(s_{t+1}| s_t,a_t)$  of moving from state $s_t\in \s_t$ to next state $s_{t+1}\in \s_{t+1}$ given action $a_t\in \A_t$.
Unless stated otherwise we consider non-randomized policies, that is $\Pi=\{\pi_1,...,\pi_T\}$ where
$\pi_t:\s_t\to \A_t$ is a (measurable) mapping form the state space to the action space.
The value $s_1$ is deterministic, initial conditions. We also use notation
$P^{s_t,a_t}_t(\cdot)= P_t(\cdot| s_t,a_t)$ for the transition kernel.

Assuming that the data process $\xi_t$ is stagewise independent, the SOC can be formulated in the MDP framework by defining transition probabilities
\begin{equation}\label{transoc}
 P^{s_t,a_t}_t(s_{t+1}\in A):=P_{t}\big (\Phi_t(s_t,a_t,\xi_t)\in A\big),\;A\in \F_{t+1},
\end{equation}
where $\xi_t\sim P_t$.
 However, there is an essential difference between the SOC and MDP modeling. In the MDP the probability law of the process is defined by the transition kernels and there is no explicitly defined random data process. The basic assumption of the SOC, used in section \ref{sec-soc},     that the   distribution of the random data process
 does not   depend on states and controls, is not directly applicable in the  MDP setting.

 The   dynamic equations for problem \eqref{mdprob-1} are:
  $V_{T+1}(s_{T+1})=  c_{T+1}(s_{T+1})$ and for $t=T,...,1$,  and $s_t\in \s_t$,
\begin{equation}\label{mdprob-2}
V_t(s_t)=\inf_{a_t\in \A_t}  \bbe^{P_t^{s_t ,a_t}}\big[
c_t(s_t, a_t ,s_{t+1})+V_{t+1}(s_{t+1}) \big],
\end{equation}
where  $\bbe^{P_t^{s_t ,a_t}}$ denotes the expectation with respect to the probability measure $P^{s_t,a_t}_t$ on $(\s_{t+1},\F_{t+1})$.
The optimal policy is given by
\begin{equation}\label{mdprob-3}
\pi_t(s_t)\in\argmin_{a_t\in \A_t}  \bbe^{P_t^{s_t ,a_t}}\big[
c_t(s_t, a_t ,s_{t+1})+V_{t+1}(s_{t+1}) \big],
\end{equation}
$t=1,...,T$.
\\

As before  we assume that with every probability measure $P_t$ on $(\s_t,\F_t)$ is associated functional $\R^{P_t}:\Z_t\to \bbr$ defined on a linear space $\Z_t$ of measurable functions $Z_t:\s_t\to \bbr$. Unless stated otherwise  we assume that   $\R^{P_t}$ satisfies the axioms of monotonicity and translation equivariance.
The counterpart of dynamic equations \eqref{mdprob-2} (compare with dynamic equations \eqref{dynrobust} for the risk averse SOC) is
\begin{equation}\label{mdprob-4}
V_t(s_t)=\inf_{a_t\in \A_t}  \R^{P_t^{s_t ,a_t}}\big[
c_t(s_t, a_t ,s_{t+1})+V_{t+1}(s_{t+1}) \big].
\end{equation}
These dynamic equations correspond to the nested formulation of the respective risk averse MDP. That is, for a policy $\pi\in \Pi$, with $a_t=\pi_t(s_t)$,   the kernels $P_t^{s_t,a_t}$ define the respective  probability distribution on the sequences $(s_1,...,s_{T+1})\in \s_1\times \cdots\times \s_{T+1}$ (Ionescu Tulcea Theorem, \cite{IT1949}). Consequently the nested functional $\cR^\pi=\R^{P_1^{s_1,\pi_1(s_1)}}\circ\cdots \circ \R^{P_{T}^{s_{T},\pi_{T}(s_{T})}}$ is defined on the history of the process $s_1,...,s_{T+1}$  (compare with \eqref{iterat-3}). That is,
for  a measurable $Z :\s_1\times\cdots\times \s_{T+1}\to \bbr$ define iteratively  $Z_{T+1}=Z$ and for $t=T,...,1$,
\begin{equation}\label{mdprob-5}
Z_t(s_1,...,s_t):= \R^{P_t^{s_t,\pi_t(s_t)}}(Z_{t+1}(s_1,...,s_t,\cdot)).
\end{equation}
Eventually for $t=1$, the value $Z_1(s_1)$ is deterministic, and
  $\cR^\pi(Z)=Z_1(s_1)$.

The  nested counterpart of problem \eqref{mdprob-1} is
\begin{equation}\label{mdprob-6}
 \min_{\pi\in \Pi} \cR^{\pi}
 \left [ \tsum_{t=1}^{T}
c_t(s_t,a_t,s_{t+1})+c_{T+1}(s_{T+1})
\right],
\end{equation}
  with  the corresponding dynamic equations given in \eqref{mdprob-4}. Let us emphasize that here the nested functional $\cR^{\pi}$  is  associated with policy $\pi\in \Pi$, and  is defined on the history of the process $s_1,...,s_{T+1}$
 with probability of moving from $s_t\in \s_t$ to $s_{t+1}\in \s_{t+1}$ given by the transition kernel  $P_t^{s_t,\pi_t(s_t)}$.
 For {\em coherent risk} measures $\R^{P_t}$, in particular for the Average Value-at-Risk,  the   nested formulation \eqref{mdprob-6}, and the respective dynamic equations \eqref{mdprob-4},  is  equivalent to the construction introduced in Ruszczy\'nski \cite{rusz2010} where it was developed in terms of the dual representation of coherent risk measures.
 In case $\R^{P_t}=\bbe^{P_t}$, problem \eqref{mdprob-6} coincides with the risk neutral problem \eqref{mdprob-1}.

 We can also consider robust   risk averse  functionals (compare with \eqref{robcoun-d}). That is, let $\cM_t$ be a set of {conditional probability measures} $P_t^{s_t,a_t}$. Then consider the following (robust) extension of dynamic equations \eqref{mdprob-4}:
  \begin{equation}\label{robdynmdp-1}
V_t(s_t)=\inf_{a_t\in \A_t} \sup_{P_t^{s_t,a_t}\in \cM_t}  \R^{P_t^{s_t ,a_t}}\big[
c_t(s_t, a_t ,s_{t+1})+V_{t+1}(s_{t+1}) \big].
\end{equation}
We can give two equivalent formulations of the corresponding robust risk averse problem. We can define the following   robust  counterpart of the  risk averse functionals
  \begin{equation}\label{robdynmdp-2}
\bfR_t^{s_t,a_t}(\cdot):=\sup_{P_t^{s_t,a_t}\in \cM_t}  \R^{P_t^{s_t ,a_t}}(\cdot),\;t=1,...,T,
\end{equation}
and the respected nested functional ${\bf R}^\pi=\bfR^{s_1,\pi_1(s_1)}\circ\cdots \circ \bfR^{s_{T},\pi_{T}(s_{T})}$
defined iteratively as in \eqref{mdprob-5}. Then the corresponding robust nested formulation is obtained by replacing $\cR^\pi$ in \eqref{mdprob-6} with ${\bf R}^\pi$, that is
\begin{equation}\label{robdynmdp-3}
 \min_{\pi\in \Pi} {\bf R}^\pi
 \left [ \tsum_{t=1}^{T}
c_t(s_t,a_t,s_{t+1})+c_{T+1}(s_{T+1})
\right].
\end{equation}

Alternatively we can formulate it as a game between the controller (decision maker) and the nature   (cf., \cite{LiSha}).
 That is, consider   the min-max problem
\begin{equation}\label{nestpr-2}
 \min_{\pi\in \Pi}\sup_{\gamma\in \Gamma} \cR^{\pi,\gamma}
 \left [ \tsum_{t=1}^{T}
c_t(s_t,\pi_t(s_t),s_{t+1})+c_{T+1}(s_{T+1})
\right],
\end{equation}
 where $\gamma_t(s_t)\in \cM_t$, $t=1,...,T$, defines the respective  policy of the nature,  and $\cR^{\pi,\gamma}$ is the nested functional  defined iteratively as in \eqref{mdprob-5} with $P_t^{s_t,a_t}=\gamma(s_t)$.

\begin{remark}
Similar to the  static formulation\footnote{In some publications such static formulations are  referred to as {\em robust} MDPs.} of  distributionally robust counterpart  of  risk-neutral MDP \cite{Nilim2005, iyen}, in the considered risk-averse setting  it is also possible to consider  the  following   counterpart of \eqref{nestpr-2},
\begin{align}\label{nestpr-static}
 \min_{\pi\in \Pi}\sup_{P_t \in \PP_t, 1 \leq t \leq T} \cR^{\pi, \{P_t\}_{t=1}^T}
 \left [ \tsum_{t=1}^{T}
c_t(s_t,\pi_t(s_t),s_{t+1})+c_{T+1}(s_{T+1})
\right],
\end{align}
where  $\cR^{\pi, \{P_t\}_{t=1}^T}(\cdot) $ refers to the nested risk functional \eqref{mdprob-6} with  fixed reference kernels $\{P_t\}_{t=1}^T$.
It should be noted that in \eqref{nestpr-static},  $\{P_t \}_{t=1}^T$ are chosen before the realization of the data process, and consequently dynamic equations \eqref{robdynmdp-1} do not apply to \eqref{nestpr-static} in general.
For dynamic equations to hold, there is a  need for additional rectangularity structure of the ambiguity sets $\PP_t$  (cf., \cite{Nilim2005, iyen,kuhn2013}), in which  case \eqref{nestpr-2}  and its static counterpart \eqref{nestpr-static} become equivalent.
\end{remark}

\subsection{Infinite horizon setting}
\label{sec-infmdp}
In the infinite horizon setting, for the counterpart of the robust risk averse problem \eqref{robdynmdp-3},  the corresponding Bellman equation can be written as
\begin{equation}\label{bellm}
V(s)=\inf_{a\in \A}\sup_{P^{s,a}\in \cM}  \R^{P^{s,a}}   \big[
c(s, a ,\cdot)+\beta V (\cdot) \big].
\end{equation}
Here $\beta\in (0,1)$ is the discount factor, $\A$ is the set of actions, $c:\s\times\A\times \s\to \bbr$ is the cost function,
$\cM$ is the  set of transition kernels, i.e., elements of   $\cM$  are conditional  probability measures  $P^{s,a}$    on the state space $(\s,\F)$, and $\R^P$ is the  risk functional. As before we assume that $\s$ and $\A$ are Polish spaces.

For bounded cost function,    equation \eqref{sec-infmdp} has a unique solution. Proof of that is rather standard (and similar to the proof for SOC discussed in section \ref{sec-infsoc}). Consider the space $\bbb$ of bounded measurable functions $g:\s\to \bbr$, equipped with the sup-norm, and Bellman operator $T:\bbb\to\bbb$ ,
 \[
T(g)(\cdot):= \inf_{a\in \A}\sup_{P^{s,a}\in \cM}  \R^{P^{s,a}}   \big[
c(s, a ,\cdot)+\beta g (\cdot) \big].
\]
  The operator $T$ has the following properties. It is monotone, i.e., if $g,g'\in \bbb$ and $g(\cdot)\ge g'(\cdot)$, then $T(g)(\cdot)\ge T(g')(\cdot)$. This follows from the monotonicity property of $\R^P$. Also   from the translation equivariance property of  $\R^P$ follows   the  property of constant shift: for any $g\in \bbb$ and $c\in \bbr$,
  $T(g+c)=T(g)+\beta c$. It follows that $T$ is a contraction mapping, i.e.,
  \[
  \|T(g)-T(g')\|_\infty \le \beta \|g-g'\|_\infty,\;g,g'\in \bbb,
\]
and hence by the Banach Fixed Point Theorem has a unique fixed point.

\bibliographystyle{plain}
\bibliography{references}

\end{document}